\documentclass{amsart}
\usepackage{amssymb,amsfonts,amsmath}
\usepackage{color,graphics}
\usepackage[matrix,arrow]{xy}

\numberwithin{equation}{section}

\theoremstyle{plain}

\newtheorem{thm}{Theorem}[section]
\newtheorem{theorem}[thm]{Theorem}

\newtheorem{corollary}[thm]{Corollary}
\newtheorem{lemma}[thm]{Lemma}

\newtheorem{proposition}[thm]{Proposition}

\newtheorem{definition}[thm]{Definition}

\newtheorem{remark}[thm]{Remark}

\newcommand{\R}{\mathbb{R}}
\newcommand{\Z}{\mathbb{Z}}
\newcommand{\C}{\mathbb{C}}

\newcommand{\half}{{\textstyle\frac{1}{2}}}

\newcommand{\htp}{\simeq}
\newcommand{\iso}{\cong}
\newcommand{\smooth}{C^\infty}

\newcommand{\A}{\mathcal{A}}
\newcommand{\B}{\mathcal{B}}

\newcommand{\tw}{\mathit{tw}}
\renewcommand{\hom}{\mathit{hom}}
\newcommand{\Hom}{\mathit{Hom}}

\begin{document}
\title[Exotic symplectic structures]{Lefschetz fibrations and exotic symplectic structures on cotangent bundles of spheres}
\author{Maksim Maydanskiy, Paul Seidel}
\maketitle

\section{Introduction}

Lefschetz fibrations provide one of the available methods for constructing symplectic structures. This paper builds on the model of \cite{maydanskiy09}, where that method was used to find a non-standard symplectic structure on the manifold obtained by attaching an $n$-handle to the cotangent bundle of the $(n+1)$-sphere (for any even $n \geq 2$). Here, we explore a somewhat more speculative idea first proposed in \cite[Section 2]{seidel07}, namely that random choices of vanishing cycles almost always lead to non-standard symplectic structures.

Generally, starting from any $2n$-dimensional Liouville manifold $M$ and an ordered collection of Lagrangian spheres $(V_1,\dots,V_r)$ in it, one constructs a $(2n+2)$-dimensional Liouville manifold $E$, namely the total space of the Lefschetz fibration with fibre $M$ and having $(V_1,\dots,V_r)$ as a basis of vanishing cycles.
In our case, we take the fibre $M = M_m$ to be the $2n$-dimensional type $(A_m)$ Milnor fibre (for some $m,n \geq 2$), whose symplectic geometry is well-studied \cite{khovanov-seidel98}. In particular, there is a map associating a Lagrangian sphere $S_{\delta} \subset M_m$ to any suitable path $\delta$ in the $(m+1)$-punctured plane. We choose our $r = m+1$ vanishing cycles $V_k = S_{\delta_k}$ as follows. The first $m$ of them are fixed, and come from an $(A_m)$ chain of paths $(\delta_1,\dots,\delta_m)$. Given those, we then allow an arbitrary $\delta_{m+1}$ (see Figure \ref{fig:random} for an example). Only the isotopy class of $\delta_{m+1}$ really matters, but still, for any fixed $m$ there is an infinite number of possible choices; they correspond bijectively to elements in the braid group $\mathit{Br}_{m+1}$ conjugate to one of the generators in the standard presentation. The homotopy type of the resulting total space $E$ is always that of $S^{n+1}$, and in fact we can be a little more precise concerning its topology:

\begin{figure}
\begin{centering}
\begin{picture}(0,0)%
\includegraphics{random.pstex}%
\end{picture}%
\setlength{\unitlength}{3355sp}%
\begingroup\makeatletter\ifx\SetFigFont\undefined%
\gdef\SetFigFont#1#2#3#4#5{%
  \reset@font\fontsize{#1}{#2pt}%
  \fontfamily{#3}\fontseries{#4}\fontshape{#5}%
  \selectfont}%
\fi\endgroup%
\begin{picture}(2363,2664)(1286,-1991)
\put(2851,-886){\makebox(0,0)[lb]{\smash{{\SetFigFont{10}{12.0}{\rmdefault}{\mddefault}{\updefault}{$\delta_5$}%
}}}}
\put(3076,164){\makebox(0,0)[lb]{\smash{{\SetFigFont{10}{12.0}{\rmdefault}{\mddefault}{\updefault}{$\delta_1$}%
}}}}
\put(1501,-61){\makebox(0,0)[lb]{\smash{{\SetFigFont{10}{12.0}{\rmdefault}{\mddefault}{\updefault}{$\delta_2$}%
}}}}
\put(1301,-586){\makebox(0,0)[lb]{\smash{{\SetFigFont{10}{12.0}{\rmdefault}{\mddefault}{\updefault}{$\delta_3$}%
}}}}
\put(2551,-1486){\makebox(0,0)[lb]{\smash{{\SetFigFont{10}{12.0}{\rmdefault}{\mddefault}{\updefault}{$\delta_4$}%
}}}}
\end{picture}%
\caption{\label{fig:random}}
\end{centering}
\end{figure}

\begin{lemma} \label{th:isotopy}
If $n$ is even, any choice of $V_{m+1}$ leads to a manifold $E$ which is diffeomorphic to the cotangent bundle $T^*\! S^{n+1}$. Moreover, under this diffeomorphism, the homotopy class of its almost complex structure turns into the standard one.
\end{lemma}

\begin{lemma} \label{th:haefliger}
Take $n$ odd (and by our previous assumption, $>1$). Choose orientations of $V_1,\dots,V_m$ such that $V_i \cdot V_{i+1} = (-1)^{\frac{1}{2} n(n+1)+1}$ for all $i<m$. Suppose that one can orient $V_{m+1}$ such that
\begin{equation} \label{eq:homology-class}
[V_{m+1}] = \sum_{i=k}^{l-1} [V_i] \in H_n(M_m) \iso \Z^m
\end{equation}
for some $k<l$. Then $E$ is diffeomorphic to $T^*\! S^{n+1}$. Moreover, under this diffeomorphism, the homotopy class of its almost complex structure turns into the standard one.
\end{lemma}

The first statement is elementary, while the second one relies on some computations in classical homotopy theory. If \eqref{eq:homology-class} is violated, then $E$ is distinguished from $T^*\!S^{n+1}$ by the intersection pairing on $H_{n+1}(E) \iso \Z$. Even though we do not consider it in the body of the paper, the construction could also be carried out for $n = 1$, leading to four-manifolds $E$ which are double branched covers of $\R^4$. However, Lemma \ref{th:haefliger} fails to hold in that dimension, because of the additional obstruction given by the fundamental group at infinity.

\begin{theorem} \label{th:main}
Among all possible choice of isotopy classes of $\delta_{m+1}$, there are exactly $\frac{1}{2} m(m+1)$ which lead to $E$ being Liouville isomorphic to $T^*\!S^{n+1}$. In all other cases, $E$ does not contain a Lagrangian sphere representing a nonzero class in $H_{n+1}(E)$, hence is not symplectomorphic to $T^*\!S^{n+1}$.
\end{theorem}

\begin{figure}
\begin{centering}
\begin{picture}(0,0)%
\includegraphics{simple.pstex}%
\end{picture}%
\setlength{\unitlength}{3355sp}%
\begingroup\makeatletter\ifx\SetFigFont\undefined%
\gdef\SetFigFont#1#2#3#4#5{%
  \reset@font\fontsize{#1}{#2pt}%
  \fontfamily{#3}\fontseries{#4}\fontshape{#5}%
  \selectfont}%
\fi\endgroup%
\begin{picture}(4563,2719)(1636,-4409)
\put(4651,-3136){\makebox(0,0)[lb]{\smash{{\SetFigFont{10}{12.0}{\rmdefault}{\mddefault}{\updefault}{$\delta_2$}%
}}}}
\put(5101,-2761){\makebox(0,0)[lb]{\smash{{\SetFigFont{10}{12.0}{\rmdefault}{\mddefault}{\updefault}{$\delta_1$}%
}}}}
\put(2251,-2761){\makebox(0,0)[lb]{\smash{{\SetFigFont{10}{12.0}{\rmdefault}{\mddefault}{\updefault}{$\delta_1$}%
}}}}
\put(2476,-3736){\makebox(0,0)[lb]{\smash{{\SetFigFont{10}{12.0}{\rmdefault}{\mddefault}{\updefault}{$\delta_3$}%
}}}}
\put(4726,-1861){\makebox(0,0)[lb]{\smash{{\SetFigFont{10}{12.0}{\rmdefault}{\mddefault}{\updefault}{$\delta_3$}%
}}}}
\put(1651,-4336){\makebox(0,0)[lb]{\smash{{\SetFigFont{10}{12.0}{\rmdefault}{\mddefault}{\updefault}{standard ($T^*\!S^{n+1}$)}%
}}}}
\put(4501,-4336){\makebox(0,0)[lb]{\smash{{\SetFigFont{10}{12.0}{\rmdefault}{\mddefault}{\updefault}{nonstandard}%
}}}}
\put(1801,-3136){\makebox(0,0)[lb]{\smash{{\SetFigFont{10}{12.0}{\rmdefault}{\mddefault}{\updefault}{$\delta_2$}%
}}}}
\end{picture}%
\caption{\label{fig:simple}}
\end{centering}
\end{figure}

Figure \ref{fig:simple} shows two choices of paths, of which the left hand one yields $T^*\!S^{n+1}$, while the right hand one corresponds to a nonstandard structure. The obvious question is whether the nonstandard structures constructed in this way are all the same, or whether there are actually infinitely many different ones among them (and if so, how they depend on the choice of path). Unfortunately, the invariant used in this paper provides no help in answering that question.

The structure of the paper is as follows. Sections \ref{sec:lefschetz}--\ref{sec:filtration} collect some properties of Floer cohomology in the context of Lefschetz fibrations. None of that is really new, hence proofs will usually be only outlined. Section \ref{sec:vanishing} discusses the main Floer cohomology computation, and Section \ref{sec:algebra} its algebraic implications. Up to this point, everything is fairly general. Then, Section \ref{sec:milnor} reviews the manifolds $M_m$, and Section \ref{sec:quiver} some elementary facts about $(A_m)$ quiver representations. As we'll see in Section \ref{sec:main}, an application of the previously introduced general ideas to this specific situation quickly leads to the desired conclusion.

{\em Acknowledgments.} The authors are indebted to Denis Auroux for valuable assistance. Remark \ref{th:wr-wr} follows suggestions by Mohammed Abouzaid. We also thank the referee for suggestions which improved the exposition. The second author was partially supported by NSF grant DMS-0652620.

\section{Lefschetz fibrations\label{sec:lefschetz}}

This section reviews some basic symplectic geometry notions used in the paper. A few references covering similar material are \cite{eliashberg-gromov93,eliashberg94,seidel01,seidel04} (this is by no means an exhaustive list).

\begin{definition} \label{def:liouville}
A Liouville manifold is a $2n$-manifold $M$ with a one-form $\theta_M$ such that $d\theta_M = \omega_M$ is symplectic, and satisfying the following additional property. Let $Z_M$ be the Liouville vector field dual to $\theta_M$. Then there is a relatively compact open subset $\Omega_M \subset M$ with smooth boundary, such that $Z_M$ points outwards along $\partial \Omega_M$, and such that the positive time flow of $Z_M$ provides a diffeomorphism $\R^+ \times \partial\Omega_M \rightarrow \overline{M \setminus \Omega_M}$ (here and below, $\R^+ = [0,\infty)$).
\end{definition}

This should be more accurately called a complete finite type Liouville manifold, but since this is the only type considered here, we omit the adjectives.

\begin{definition}
Let $M$ and $N$ be Liouville manifolds. A Liouville isomorphism is a diffeomorphism $\psi: M \rightarrow N$ such that $\psi^*\theta_N - \theta_M$ is the derivative of a compactly supported function.
\end{definition}

We remind the reader of some basic implications of these definitions. If $\Omega_M \subset M$ is as before, then $\theta_M|\partial \Omega_M$ is a contact one-form, and any two hypersurfaces obtained in this way are canonically contactomorphic (by going along the flow lines of $Z_M$). Hence, $\partial \Omega_M$ is sometimes called the boundary at infinity of $M$. Liouville isomorphisms induce contact isomorphisms between the respective boundaries at infinity.

Our next task is to introduce Lefschetz fibrations, in a version which is suitable for our constructions. Since both the fibres and the base will be noncompact, some control on the geometry near infinity has to be imposed, and we'll set that up first. Fix a $2n$-dimensional Liouville manifold $(M,\omega_M = d\theta_M)$ and a compactly supported Liouville automorphism $\mu$ of $M$. Take an open subset with smooth boundary in $M$, whose complement is compact, and let $U$ be the closure of that subset. We can choose $U$ such that $\mu|U$ is the identity, and so that the function $k_\mu$ satisfying $\mu^*\theta_M = \theta_M + dk_\mu$ vanishes on $U$. Let $T_\mu = (\R \times M)/(t,x) \sim (t-2\pi, \mu(x))$ be the mapping torus of $\mu$. Define
\begin{equation} \label{eq:fake-pi}
\pi_{\tilde{E}}: \tilde{E} = (\C \times {U}) \cup_{\R^+ \times S^1 \times {U}} (\R^+ \times T_\mu)  \longrightarrow \C
\end{equation}
as follows. The identification of the two pieces takes $(\exp(s+it),x)$ to $(s,t,x)$ (technically, this makes $\tilde{E}$ into a manifold with ``concave'' codimension $2$ corners). The map is $\pi_{\tilde{E}}(z,x) = z$ on $\C \times {U}$, and $\pi_{\tilde{E}}(s,t,x) = \exp(s+it)$ on $\R^+ \times T_\mu$. $\tilde{E}$ carries a symplectic form $\omega_{\tilde{E}}$, which is equal to the standard product form $d\mathrm{re}(z) \wedge d\mathrm{im}(z) + \omega_M$ on $\C \times {U}$, and to $e^{2s} ds \wedge dt + \omega_M$ on $\R^+ \times T_\mu$. This has a one-form primitive $\theta_{\tilde{E}}$, which is equal to $\half(\mathrm{re}(z)\, d\mathrm{im}(z) - \mathrm{im}(z)\, d\mathrm{re}(z)) + \theta_M$ on $\C \times {U}$, and to $\frac{1}{2} e^{2s} dt + \theta_M + d( (t/2\pi) k_\mu)$ on $\R^+ \times T_\mu$.

\begin{definition} \label{th:fibration}
A Lefschetz fibration (with fibre $M$ and outer monodromy $\mu$) is a smooth map $\pi_E: E \rightarrow \C$, where $E$ is a $(2n+2)$-manifold together with an exact symplectic form $\omega_E = d\theta_E$, satisfying the following properties. At each regular point, $\mathrm{ker}(D\pi_E) \subset TE$ is a symplectic subspace. Besides that we have the Lefschetz condition, which says that locally near each critical point, our fibration is modelled after the complex function
\begin{equation} \label{eq:morse-chart}
\C^{n+1} \longrightarrow \C, \quad
x \longmapsto x_1^2 + \cdots + x_{n+1}^2 + \mathit{constant},
\end{equation}
with the symplectic form being equal to the standard constant K{\"a}hler form in the same coordinates. For simplicity, we also ask that each fibre contain at most one critical point.

Finally, the geometry near infinity is controlled by the following requirement. There is an open subset of $E$, whose closure is a compact manifold with corners, and whose complement can be identified with $\tilde{E}$. This identification should be such that the restrictions of $(\pi_E,\omega_E,\theta_E)$ equal $(\pi_{\tilde{E}},\omega_{\tilde{E}},\theta_{\tilde{E}})$.
\end{definition}

Again, these should be called exact symplectic (or Liouville) Lefschetz fibrations, but we omit all adjectives since that is the only class of Lefschetz fibrations relevant to this paper. If we look at the model \eqref{eq:fake-pi}, the Liouville vector field dual to $\theta_{\tilde{E}}$ is
\begin{equation}
Z_{\tilde{E}} = \begin{cases} \half z \partial_z + Z_M & \text{on $\C \times {U}$,} \\
(\half + e^{-2s} \textstyle\frac{k_\mu}{2\pi}) \partial_s + Z_M - \textstyle\frac{t}{2\pi} X_{k_\mu} &
\text{on $[0,\infty) \times T_\mu$,} \end{cases}
\end{equation}
where $X_{k_\mu}$ is the Hamiltonian vector field of that function. $Z_{\tilde{E}}$ points strictly outwards along the hypersurface $|\pi_{\tilde{E}}(x)| \geq r$ for sufficiently large $r$, and the same along $\C \times \partial \Omega_M$, where $\Omega_M \subset M$ is as in Definition \ref{def:liouville} and sufficiently large (so that $\partial \Omega_M \subset {U}$). By taking the subset of $E$ bounded by the combination of these hypersurfaces, and rounding off the corners, one gets a relatively compact open subset $\Omega_E$ such that $Z_E$ points outwards along its boundary. It is not difficult to show that every point outside $\Omega_E$ will be carried into that set by the flow of $Z_E$ for some negative time. As a consequence, the total space $E$ is a Liouville manifold.

\begin{definition}
A vanishing path for a Lefschetz fibration $\pi_E: E \rightarrow \C$ is a properly embedded path $\beta: \R^+ \rightarrow \C$ such that $\beta(0)$ is a singular value of $\pi_E$, and all the other $\beta(t)$ are regular values. Moreover,
\begin{equation} \label{eq:straight}
\beta(t) = t a \text{ for all $t \gg 0$}.
\end{equation}
Here $a = a_\beta \in S^1 \subset \C$ is some angle, which can be different for different vanishing paths.
\end{definition}

To any vanishing path $\beta$ one associates its Lefschetz thimble $\Delta_\beta$, which is a properly embedded Lagrangian submanifold of $E$ with $\pi_E(\Delta_\beta) = \beta(\R^+)$ (see e.g.\ \cite[Lemma 1.13]{seidel01} for the exact definition). Now suppose that we are given a basis of vanishing paths $(\gamma_1,\dots,\gamma_r)$. This means that $r$ equals the number of critical values, that all $\gamma_k$ go to infinity in positive real direction, and that they satisfy certain intersection conditions (see Figure \ref{fig:basis2} for an illustration; the notion is a classical one in singularity theory, appearing for instance as ``distinguished system'' in \cite[p.\ 60]{aglv}). Take the fibre $E_z$ for $z$ real and $\gg 0$, and identify it with $M$ in the canonical way inherited from $\tilde{E}$. The intersections of the Lefschetz thimbles $\Delta_{\gamma_k}$ with that fibre give rise to a collection of vanishing cycles $(V_1,\dots,V_r)$ in $M$, which is again called a basis. These are Lagrangian spheres in the following sense:

\begin{figure}
\begin{centering}
\begin{picture}(0,0)%
\includegraphics{basis2.pstex}%
\end{picture}%
\setlength{\unitlength}{3947sp}%
\begingroup\makeatletter\ifx\SetFigFont\undefined%
\gdef\SetFigFont#1#2#3#4#5{%
  \reset@font\fontsize{#1}{#2pt}%
  \fontfamily{#3}\fontseries{#4}\fontshape{#5}%
  \selectfont}%
\fi\endgroup%
\begin{picture}(4232,2086)(444,-1925)
\put(1126,-511){\makebox(0,0)[lb]{\smash{{\SetFigFont{11}{13.2}{\rmdefault}{\mddefault}{\updefault}{$\gamma_3$}%
}}}}
\put(2251,-961){\makebox(0,0)[lb]{\smash{{\SetFigFont{11}{13.2}{\rmdefault}{\mddefault}{\updefault}{$\gamma_1$}%
}}}}
\put(1951, 14){\makebox(0,0)[lb]{\smash{{\SetFigFont{11}{13.2}{\rmdefault}{\mddefault}{\updefault}{$\gamma_4$}%
}}}}
\put(1426,-1861){\makebox(0,0)[lb]{\smash{{\SetFigFont{11}{13.2}{\rmdefault}{\mddefault}{\updefault}{$\gamma_2$}%
}}}}
\end{picture}%
\caption{\label{fig:basis2}}
\end{centering}
\end{figure}

\begin{definition}
By a Lagrangian sphere in $M$, we mean a submanifold $V \subset M$ such that $\theta_M|V$ is exact, together with a diffeomorphism $S^n \rightarrow V$, the latter fixed up to isotopy and composition with elements of $O(n+1)$.
\end{definition}

In \cite{seidel04} these were called framed exact Lagrangian spheres, but all spheres considered in this paper naturally come with this structure, so we feel justified in shortening the terminology. The important converse to this observation is that, given $M$ and an arbitrary ordered collection of Lagrangian spheres $(V_1,\dots,V_r)$, one can construct a Lefschetz fibration $\pi_E: E \rightarrow \C$ with fibre $M$ and for which $(V_1,\dots,V_r)$ is a basis of vanishing cycles.
We refer to \cite[Section 16e]{seidel04} for a sketch of the construction, and to \cite[Proposition 1.11]{seidel01} for a more in-depth discussion of the basic case of a single vanishing cycle.

\begin{remark}
Equivalently, one can think of $E$ as being obtained by taking a compact subset of $M \times \C$, attaching Weinstein handles \cite{weinstein91} to $r$ Legendrian spheres in its boundary, then making the result non-compact again by attaching an infinite cone.
\end{remark}

If we fix $M$ and the isotopy classes of the vanishing cycles, the resulting Lefschetz fibration $\pi_E: E \rightarrow \C$ is unique up to deformation in an appropriate sense, which in particular implies that the total space $E$ is unique up to Liouville isomorphism. Moreover, there is an action of the braid group $\mathit{Br}_r$ on the set of bases of vanishing cycles by Hurwitz moves, which also leaves the isomorphism class of the total space invariant. See for instance \cite[Section 16d]{seidel04}.

Given a Lefschetz fibration $\pi_E: E \rightarrow \C$, consider a function of the form $H(y) = \psi(\half |\pi_E(y)|^2)$, where
\begin{equation} \label{eq:psi}
\text{$\psi(r) = 0$ for $r \leq \half$, and $\psi'(r) = 1$ for $r \gg 0$.}
\end{equation}
We denote by $(\Phi^\alpha)$ the flow of $H$, which is well-defined for all times $\alpha \in \R$. This flow is nontrivial only on the part of $E$ corresponding to $\R^+ \times S^1 \times M \subset \tilde{E}$, where it is given by $(s,t,x) \mapsto (s,t+\alpha \psi'(\half e^{2s}),x)$. It is fibered over a flow $(\phi^\alpha)$ on the base, which is the Hamiltonian flow of the function $h(z) = \psi(\half|z|^2)$ with respect to the standard constant symplectic form. Finally, if we take $\Phi^{2\pi}$ and restrict it to fibres $E_z \iso M$ with $|z| \gg 0$, it is fibre-preserving and fibrewise equal to the monodromy $\mu$. Given any vanishing path $\beta$, we write $\beta^\alpha = \phi^\alpha \circ \beta$. Similarly, for any Lagrangian submanifold $L \subset E$, let
\begin{equation} \label{eq:isotopy}
L^\alpha = \Phi^\alpha(L).
\end{equation}
If $L = \Delta_{\beta}$, then $L^\alpha = \Delta_{\beta^\alpha}$ is again a Lefschetz thimble.

\section{Floer cohomology\label{sec:wrapped}}

This section outlines the structure of Lagrangian Floer cohomology in the context of Lefschetz fibrations. Since we impose strong exactness conditions, this is technically rather undemanding, and we will give details only when they are particularly relevant to the intended application (see also \cite[Section 6]{maydanskiy09} for a closely related exposition).

We want to consider two classes of exact Lagrangian submanifolds $L \subset E$: closed ones and Lefschetz thimbles. Take two such submanifolds $(L_0,L_1)$, and if they are both Lefschetz thimbles $L_k = \Delta_{\beta_k}$, assume that the angles $a_k = a_{\beta_k}$ from \eqref{eq:straight} are different:
\begin{equation} \label{eq:unequal-angle}
a_0 \neq a_1.
\end{equation}
Then their Floer cohomology $HF^*(L_0,L_1)$ is well-defined. From now on, we also make the standard assumption that the anticanonical bundle $K_E^{-1} = \lambda^{n+1}_\C(TE)$ should be trivial, and that all Lagrangian submanifolds involved come with gradings \cite{seidel99}, which result in $HF^*(L_0,L_1)$ being $\Z$-graded. As for coefficients, all our Floer cohomology groups will be defined over $\Z/2$, to avoid sign considerations (this is particularly important later when we quote results from \cite{khovanov-seidel98}, where the sign issues have not been explored).

If $L_0^\alpha$ is as in \eqref{eq:isotopy}, the Floer cohomology groups $HF^*(L_0^\alpha,L_1)$ are defined for all $\alpha$ such that $e^{i\alpha}a_0 \neq a_1$ (where this condition is understood to be vacuous unless both submanifolds are Lefschetz thimbles). Moreover, there are canonical continuation maps
\begin{equation} \label{eq:continuation}
HF^*(L_0^{\alpha_-},L_1) \longrightarrow HF^*(L_0^{\alpha_+},L_1)
\end{equation}
for all values $\alpha_- \leq \alpha_+$ such that both sides make sense, and these form a directed system. We define the wrapped Floer cohomology of $(L_0,L_1)$ to be the direct limit
\begin{equation} \label{eq:direct-limit}
HW^*(L_0,L_1) = \underrightarrow{\mathit{lim}}_{\alpha}\, HF^*(L_0^\alpha,L_1).
\end{equation}

\begin{remark} \label{th:wr-wr}
The terminology ``wrapped Floer cohomology'' was introduced in \cite{fukaya-seidel-smith07b,abouzaid-seidel07} in a slightly different context; namely, for exact Lagrangian submanifolds in Liouville manifolds which are properly embedded, and are tangent to the Liouville flow outside a compact subset, hence give rise to a Legendrian submanifold of the boundary at infinity. To explain the relation between the two notions, consider a single Lefschetz thimble $L = \Delta_\beta \subset E$. Suppose that the underlying vanishing path has angle $a_\beta = 0$, and that the associated vanishing cycle $V \subset E_{\beta(t)} \iso M$, for $t \gg 0$, satisfies $\theta_M|V = 0$ (the latter condition can always be arranged by modifying the given $\theta_M$ and $\theta_E$, which does not affect Floer cohomology groups). Fix $r \gg 0$, and consider the hypersurface $\{r\} \times T_\mu \subset \tilde{E}$, which is just the subset of points where $|\pi_{\tilde{E}}(s,t,x)| = e^r$. Then
\begin{equation} \label{eq:r-tilde}
R_{\tilde{E}} = (\half e^{2r} + \textstyle\frac{1}{2\pi} k_\mu)^{-1} \partial_t
\end{equation}
is a vector field tangent to the characteristic foliation, and satisfies $\theta_{\tilde{E}}(R_{\tilde{E}}) = 1$. When we construct $\Omega_E$ by rounding off corners, this can be done in such a way that $L \cap \partial \Omega_E = L \cap E_{e^r} = V$ is Legendrian, and so that the Reeb flow applied to $L \cap \partial \Omega_E$ equals the flow of \eqref{eq:r-tilde}. The wrapped Floer cohomology in the sense of \cite{fukaya-seidel-smith07b,abouzaid-seidel07} uses this Reeb flow to form a direct limit. To obtain the isomorphism between that and $HW^*(L,L)$ as defined in \eqref{eq:direct-limit}, one defines continuation maps which intertwine the groups in the two direct systems. For that, it is crucial that $D\pi_{\tilde{E}}(R_{\tilde{E}})$ is, at every point, a positive multiple of the rotational vector field on the base (this was pointed out to the authors by Abouzaid; there is a similar argument in \cite{mclean09} for symplectic cohomology, but that is considerably more complicated, since one needs to take into account the entire boundary at infinity). We do not want to pursue this further, and instead stick to  \eqref{eq:direct-limit} as the definition of wrapped Floer cohomology, which is sufficient for our purpose (essentially the same solution is adopted in \cite[Section 6]{maydanskiy09}).
\end{remark}

Here are some general properties of \eqref{eq:continuation}. If either $L_0$ or $L_1$ is compact, the directed system is constant, so $HW^*(L_0,L_1) = HF^*(L_0,L_1)$ is ordinary Floer cohomology. Furthermore, the direct limit is compatible with the product structure, which means that there are induced associative products $HW^*(L_1,L_2) \otimes HW^*(L_0,L_1) \rightarrow HW^*(L_0,L_2)$. These are unital, and their unit elements arise as follows. For $\alpha \in (0,2 \pi)$, one has a version of the Piunikhin-Salamon-Schwarz \cite{pss} isomorphism $H^*(L;\Z/2) \iso HF^*(L^\alpha,L)$, and in particular a distinguished element $1 \in H^0(L;\Z/2) \iso HF^0(L^\alpha,L)$. The unit element in $HW^*(L,L)$ is the image of this element under \eqref{eq:direct-limit}. Associativity and unitality mean that one can introduce a wrapped version of the Donaldson-Fukaya category, having Lagrangian submanifolds as objects and wrapped Floer groups as morphisms.

\begin{lemma} \label{th:intersection-number}
Suppose that for some Lefschetz thimble $L_1$, we have $HW^*(L_1,L_1) = 0$. Then, for every closed exact $L_0$, the intersection number $L_0 \cdot L_1$ is zero.
\end{lemma}

\proof This is a formal consequence of our previous remarks. Because of the product structure, vanishing of $HW^*(L_1,L_1)$ implies vanishing of $HW^*(L_0,L_1)$ for any $L_0$. But if $L_0$ is closed, the Euler characteristic of the latter group is $\pm L_0 \cdot L_1$. \qed

\section{A spectral sequence\label{sec:filtration}}

We'll now focus on Floer cohomology for Lefschetz thimbles. By definition, such thimbles project to paths in the base $\C$ of the Lefschetz fibration. If one chooses appropriate almost complex structures, the pseudo-holomorphic strips which define the Floer differential also project to holomorphic strips in $\C$. This choice of almost complex structure is not generic, of course, but the idea can nevertheless be used for some partial computations.

Consider two vanishing paths $\beta_0,\beta_1$ with the following properties. If $\beta_0(0) = \beta_1(0)$, then $\beta'_0(0)$ and $\beta'_1(0)$ should not be positive multiples of each other. Everywhere else, $\beta_0$ and $\beta_1$ should intersect transversally, which in particular implies \eqref{eq:unequal-angle}. Consider the intersection points $z \in \beta_0(\R^+) \cap \beta_1(\R^+)$. We say that $z_- > z_+$ if there is a finite sequence of non-constant holomorphic maps
\begin{equation} \label{eq:order}
\left\{\begin{aligned}
& w_1,\dots,w_r: \R \times [0,1] \rightarrow \C, \\
& \textstyle \int_{\R \times [0,1]} |dw_k|^2 < \infty, \\
& w_k(\R \times \{0\}) \subset \beta_0(\R^+), \\
& w_k(\R \times \{1\}) \subset \beta_1(\R^+), \\
& \textstyle\lim_{s \rightarrow +\infty} w_k(s,\cdot) = \lim_{s \rightarrow -\infty} w_{k+1}(s,\cdot), \\
& \textstyle\lim_{s \rightarrow -\infty} w_1(s,\cdot) = z_-, \\
& \textstyle\lim_{s \rightarrow +\infty} w_r(s,\cdot) = z_+.
\end{aligned} \right.
\end{equation}
From now on, we partition the intersections points into subsets $I_0,\dots,I_d$, such that $z_- > z_+$ implies that $z_- \in I_i$, $z_+ \in I_j$ with $i > j$. Moreover, if our paths have the same endpoint, that point will be denoted by $b$.

Write $L_0 = \Delta_{\beta_0}$, $L_1 = \Delta_{\beta_1}$. At every intersection point $z \neq b$, we have vanishing cycles $V_{z,0},V_{z,1} \subset E_z$ ($E_z \iso M$, but not canonically so), which are simply the parts of the Lefschetz thimbles lying in that fibre. We can then consider their Floer cohomology in $E_z$, denoted by $H_z = HF^*(V_{z,0},V_{z,1})$ (the gradings of the vanishing cycles are adjusted in such a way that the indices of intersection points agree with those in the total space). In the remaining case $z = b$, $L_0$ and $L_1$ intersect transversally at the single singular point of $E_{z}$, and we set $H_z = \Z/2$, concentrated in the degree given by the Maslov index of that point.

\begin{proposition} \label{th:ss}
There is a spectral sequence converging to $HF^*(L_0,L_1)$, whose starting page has
\begin{equation} \label{eq:e1}
E_1^{pq} = \begin{cases}
\bigoplus_{z \in I_p} H_z^{p+q} & \text{$0 \leq p \leq d$}, \\
0 & \text{$p < 0$ or $p > d$}.
\end{cases}
\end{equation}
\end{proposition}

The proof is similar to Morse-Bott situations such as \cite{pozniak}, and technically even somewhat simpler.

\proof As a preliminary step, by slightly perturbing the symplectic connection on the regular part of $\pi_E: E \rightarrow \C$, we can achieve that the Lefschetz thimbles for our given paths intersect transversally. This perturbation is a compactly supported exact change of the symplectic form on $E$, which does not affect $HF^*(L_0,L_1)$ or $HF^*(V_{z,0},V_{z,1})$. To simplify the discussion, we will assume from now on that the original symplectic form already had this transversality property.

Let $CF^*(L_0,L_1)$ be the Floer complex, generated by intersection points $x \in L_0 \cap L_1$. Fix a family $(J_t)_{0 \leq t \leq 1}$ of almost complex structures on $E$ such that $\pi_E$ is $J_t$-holomorphic for each $t$. We say that $x_- > x_+$ if there is a finite sequence of non-constant pseudo-holomorphic strips
\begin{equation}
\left\{\begin{aligned}
& u_1,\dots,u_r: \R \times [0,1] \rightarrow E, \\
& \partial_s u_k + J_t(u_k) \partial_t u_k = 0, \\
& \textstyle \int_{\R \times [0,1]} \|du_k\|^2 < \infty, \\
& u_k(\R \times \{0\}) \subset L_0, \\
& u_k(\R \times \{1\}) \subset L_1, \\
& \textstyle\lim_{s \rightarrow +\infty} u_k(s,\cdot) = \lim_{s \rightarrow -\infty} u_{k+1}(s,\cdot), \\
& \textstyle\lim_{s \rightarrow -\infty} u_1(s,\cdot) = x_-, \\
& \textstyle\lim_{s \rightarrow +\infty} u_r(s,\cdot) = x_+.
\end{aligned}
\right.
\end{equation}
By projecting to the base, one sees that $x_- > x_+$ implies that $x_\pm$ either lie in the same fibre or else satisfy $\pi_E(x_-) > \pi_E(x_+)$.

The family $(J_t)$ is not generic, and in order to satisfy the transversality requirements in the definition of the Floer differential, one generally needs to perturb it slightly to some $(\tilde{J}_t)$. Nevertheless, Gromov compactness ensures that as long as the perturbation is sufficiently small, pseudo-holomorphic strips with endpoints $x_\pm$ can only exist if $x_- > x_+$.  As a consequence, if we consider the descending filtration $F^*$ of $CF^*(L_0,L_1)$ such that $F^p$ is generated by intersection points in $\pi_E^{-1}(I_p \cup I_{p+1} \cup \cdots \cup I_d)$, then the Floer differential preserves that filtration, and moreover, the induced differential on the graded space $F^p/F^{p+1}$ splits into a direct sum indexed by points $z \in I_p$. This automatically gives rise to a spectral sequence of the general form \eqref{eq:e1}. The remaining step is to determine the precise nature of the pieces $H_z$ which make up the cohomology of $F^p/F^{p+1}$. For $z = b$, the chain complex underlying $H_z$ has a single generator, so the differential automatically vanishes.

One can arrange the original family $(J_t)$ so that its restriction to the fibre $E_z$ gives rise to regular moduli spaces of pseudoholomorphic strips in that fibre, for any intersection point $z \neq b$. These strips are then automatically regular in the total space as well. Concretely, let $u: \R \times [0,1] \rightarrow E_z \subset E$ be such a strip. Its linearized operator, as a map to $E$, is a Fredholm operator $D_u: \mathcal{H}_1 \rightarrow \mathcal{H}_0$. On the other hand, if we consider $u$ as a map to $E_z$, its linearization is described by the restriction of $D_u$ to subspaces $\bar{\mathcal{H}}_1 \rightarrow \bar{\mathcal{H}}_0$. Via projection to the base, the quotient $\mathcal{H}_1/\bar{\mathcal{H}}_1$ can be identified with the space of $W^{1,q}$ (for some $q>2$) functions $\xi: \R \times [0,1] \rightarrow \C$ satisfying boundary conditions $\xi(\R \times \{0\}) \subset \sigma_0\R$, $\xi(\R \times \{1\}) \in {\sigma_1}\R$. Here, $\sigma_k\R \subset \C$ are the tangent spaces of $\beta_k(\R^+)$ at $z$, hence transverse by assumption. Similarly, $\mathcal{H}_0/\bar{\mathcal{H}}_0$ can be identified with the space of all $L^q$ functions $\R \times [0,1] \rightarrow \C$, and the quotient map induced by $D_u$ is the standard Cauchy-Riemann operator $\bar\partial$, which is invertible. This implies that regularity in $E_z$ and in $E$ are equivalent, as claimed.
Because of this regularity, the zero-dimensional moduli spaces of pseudo-holo\-mor\-phic strips for $(\tilde{J}_t)$ which define the differential on $F^p/F^{p+1}$ correspond bijectively to those for $(J_t)$, hence $H_z$ ($z \neq b$) is really the Floer cohomology in the fibre.
\qed

By similar means, one can show that $CF^*(L_0,L_1)$, considered as a filtered chain complex, is independent of the choice of almost complex structures up to isomorphism. Hence, the whole spectral sequence is canonical.
\begin{figure}
\begin{centering}
\begin{picture}(0,0)%
\includegraphics{wrapped-path.pstex}%
\end{picture}%
\setlength{\unitlength}{3947sp}%
\begingroup\makeatletter\ifx\SetFigFont\undefined%
\gdef\SetFigFont#1#2#3#4#5{%
  \reset@font\fontsize{#1}{#2pt}%
  \fontfamily{#3}\fontseries{#4}\fontshape{#5}%
  \selectfont}%
\fi\endgroup%
\begin{picture}(4737,3548)(-61,-2054)
\put(3676,-361){\makebox(0,0)[lb]{\smash{{\SetFigFont{11}{13.2}{\rmdefault}{\mddefault}{\updefault}{$z_2$}%
}}}}
\put(4201,-711){\makebox(0,0)[lb]{\smash{{\SetFigFont{11}{13.2}{\rmdefault}{\mddefault}{\updefault}{$\beta_1$}%
}}}}
\put(3751,164){\makebox(0,0)[lb]{\smash{{\SetFigFont{11}{13.2}{\rmdefault}{\mddefault}{\updefault}{$\beta_0$}%
}}}}
\put(3451,1139){\makebox(0,0)[lb]{\smash{{\SetFigFont{11}{13.2}{\rmdefault}{\mddefault}{\updefault}{$\beta_0^{4\pi}$}%
}}}}
\put(1575,239){\makebox(0,0)[lb]{\smash{{\SetFigFont{11}{13.2}{\rmdefault}{\mddefault}{\updefault}{$b$}%
}}}}
\put(2850,-626){\makebox(0,0)[lb]{\smash{{\SetFigFont{11}{13.2}{\rmdefault}{\mddefault}{\updefault}{$z_1$}%
}}}}
\end{picture}%

\caption{\label{fig:wrapped-path}}
\end{centering}
\end{figure}

We will need a related consideration concerning continuation maps. This requires some technical assumptions:
\begin{equation} \label{eq:psi2}
\parbox{32em}{
$\beta$ is a vanishing path, which outside the unit disc is a straight half-line. $L = \Delta_\beta$ is its Lefschetz thimble. Moreover, the function from \eqref{eq:psi} should satisfy $\psi''(r) > 0$ for all $r$ such that $\psi'(r) \neq 0,1$.}
\end{equation}
Take $\beta_1 = \beta$, and consider another path $\beta_0$ of the following kind. $\beta_0$ is isotopic to $\beta$ inside the class of vanishing paths; it meets $\beta$ only at the common endpoint $b$, where their oriented tangent directions are different; and outside the unit disc it equals $e^{i\alpha}\beta$ for some $\alpha \in (0,2\pi)$. Write $L_k = \Delta_{\beta_k}$. Fix some positive integer $d$. The wrapped path $\beta_0^{2\pi d}$ intersects $\beta_1$ at the origin and at other points $z_1,\dots,z_d$  (see Figure \ref{fig:wrapped-path}), where the intersections are automatically transverse because of \eqref{eq:psi2}. These points satisfy $b > z_1 > \cdots > z_d$ in the sense of \eqref{eq:order}. Consider the spectral sequence obtained by setting $I_d = \{b\}$, and $I_k = \{z_{d-k}\}$ for $k<d$. As part of that spectral sequence, we have an edge homomorphism
\begin{equation} \label{eq:edge}
H_b = \Z/2 \longrightarrow HF^0(L_0^{2\pi d},L_1).
\end{equation}
Recall that by construction, $HF^0(L_0^{2\pi d},L_1) \iso HF^0(L^{2\pi d + \alpha},L)$. Hence, it carries a canonical element, which is the image of the nontrivial generator of $HF^0(L_0,L_1) \iso HF^0(L^\alpha,L) = \Z/2$ under \eqref{eq:continuation}. In the limit \eqref{eq:direct-limit}, this gives rise to the unit element of $HW^*(L,L)$.

\begin{proposition} \label{th:element}
The canonical element of $HF^0(L_0^{2\pi d},L_1)$ is given by \eqref{eq:edge}.
\end{proposition}

\proof The argument hinges on the chain level realization of the continuation map,
\begin{equation} \label{eq:cont-d}
\Z/2 \iso CF^*(L_0,L_1) \longrightarrow CF^*(L_0^{2\pi d},L_1).
\end{equation}
Fix some $\eta \in \smooth(\R,\R)$ such that
\begin{equation} \label{eq:eta}
\text{$\eta(s) = 2\pi d$ for $s \ll 0$, $\eta(s) = 0$ for $s \gg 0$, and $\eta'(s) \leq 0$ everywhere.}
\end{equation}
Then, \eqref{eq:cont-d} is defined by counting solutions of
\begin{equation} \label{eq:cont-u}
\left\{
\begin{aligned}
& u: \R \times [0,1] \longrightarrow E, \\
& \textstyle\int_{\R \times [0,1]} \|\partial_t u\|^2 < \infty, \\
& u(s,0) \in L_0^{\eta(s)}, \;\; u(\R \times \{1\}) \subset L_1, \\
& \partial_s u + J_{s,t}(u) \partial_t u - (1-t)\eta'(s)X_H(u) = 0.
\end{aligned}
\right.
\end{equation}
Assume temporarily that all $J_{s,t}$ have the property that $\pi_E$ is $J_{s,t}$-holomorphic. Then, for every $u$ as in \eqref{eq:cont-u}, the projection $w = \pi_E(u)$ is a solution of
\begin{equation} \label{eq:cont-w}
\left\{
\begin{aligned}
& w: \R \times [0,1] \longrightarrow \C, \\
& \textstyle\int_{\R \times [0,1]} |\partial_t w|^2 < \infty, \\
& w(s,0) \in \beta_0^{\eta(s)}(\R^+), \;\; w(\R \times \{1\}) \subset \beta_1(\R^+), \\
& \partial_s w + i\partial_t w - (1-t)\eta'(s) X_h(w) = 0.
\end{aligned}
\right.
\end{equation}
Take a solution of \eqref{eq:cont-w}, and write it as $w(s,t) = \phi^{(1-t)\eta(s)}(v(s,t))$, so that the boundary values of $v$ lie on $\beta_0$ and $\beta_1$, respectively. If $w$ has limit $b$ as $s \rightarrow \infty$, and limit $z_k$ as $s \rightarrow -\infty$, then
\begin{equation} \label{eq:energy}
\begin{aligned}
0 & \leq \textstyle \int_{\R \times [0,1]} |\partial_t w|^2 \\ & =
\textstyle \int_{\R \times [0,1]} \omega_\C (\partial_s w - (1-t)\eta'(s)X_h(w), \partial_t w) \\ & =
\textstyle \int_{\R \times [0,1]} \omega_\C (\partial_s v, \partial_t v - \eta(s)X_h(v)) \\ & = 
\textstyle \int_{\R \times [0,1]} v^* \omega_\C - \eta(s) dh(\partial_s v) \\ & \leq
\textstyle \int_{\R \times [0,1]} v^* \omega_\C - \partial_s (\eta(s) h(v)) \\ & =
\textstyle \Big(\int_{\R \times [0,1]} v^* \omega_\C \Big) + 2\pi d \cdot \psi(\half |z_k|^2).
\end{aligned}
\end{equation}
Here, $\omega_\C$ is the standard symplectic form. The second inequality holds because $\eta'(s) h(z) = \eta'(s) \psi(\half |z|^2)$ is non-positive everywhere, by assumptions \eqref{eq:psi2} and \eqref{eq:eta}.

The rest of the computation is elementary. The first term in the last line of \eqref{eq:energy} is $A - \pi k + (\alpha/2 - \pi k) (|z_k|^2 - 1)$, where $A \in (0,\pi)$ is the area of the curvilinear triangle bounded by (parts of) $\beta_0$, $\beta_1$, and the unit circle. On the other hand,
\begin{equation}
\textstyle 2\pi d \cdot \psi(\half |z_k|^2) = \int_1^{|z_k|} 2\pi d \cdot \psi'(\half r^2) r \, dr
\end{equation}
where by \eqref{eq:psi2} $2 \pi d \cdot \psi'(\half r^2)$ is a monotone function with values going from zero to $2\pi k - \alpha$ in our domain of integration. Hence, $2 \pi d \cdot \psi(\half |z_k|^2) \leq (\pi k - \alpha/2)(|z_k|^2 - 1)$. Taking everything together, we find that
\begin{equation}
0 \leq \textstyle \int_{\R \times [0,1]} |\partial_t w|^2 \leq A - \pi k,
\end{equation}
an obvious contradiction. Hence, any solution of \eqref{eq:cont-w} necessarily has limits $b$ at both ends. The same computation shows that such solutions are also necessarily constant. This in turn means that all solutions of \eqref{eq:cont-u} are contained in the fibre $E_b$. Since $X_H = 0$ in that fibre, they are actually pseudo-holomorphic maps, hence also constant by a straightforward energy argument. Finally, we should note that the constant map is a regular solution.

This holds for our non-generic choice of $(J_{s,t})$. However, the same Gromov compactness argument as in Proposition \ref{th:ss} shows that for any sufficiently small perturbation $(\tilde{J}_{s,t})$, the continuation map equation still has the constant map at the singular point in $E_b$ as its only solution. Hence, the homomorphism defined in this way indeed agrees with \eqref{eq:edge}. \qed

Proposition \ref{th:element} has an analogue for symplectic cohomology $SH^*(E)$, which is maybe better known. Namely, if one defines $SH^*(E)$ using a suitable Hamiltonian function $K$, then the image of the canonical map $H^*(E) \rightarrow SH^*(E)$ is represented by the critical points of $K$, seen as (constant) one-periodic orbits of the flow of $X_K$. This is implicit in \cite{viterbo97a}, for instance.

\section{A vanishing criterion\label{sec:vanishing}}

This section contains our main Floer cohomology computation. It takes place in the framework of the spectral sequence from Proposition \ref{th:ss}, and in fact in a particularly simple special case where only one differential in that spectral sequence can be nontrivial. However, computing that differential requires additional geometric considerations.

Let $\pi_E: E \rightarrow \C$ be a Lefschetz fibration as before, $(\gamma_1,\dots,\gamma_r)$ a basis of vanishing paths, and $(V_1,\dots,V_r)$ the corresponding basis of vanishing cycles. To each $V_k \subset M$ we can associate its Dehn twist $\tau_{V_k}$, which is an exact symplectic automorphism of $M$ (and graded in a canonical way); the Picard-Lefschetz formula then tells us that the outer monodromy $\mu$ is isotopic to $\tau_{V_1}\cdots\tau_{V_r}$. From the long exact sequence \cite{seidel01} we know that there are preferred elements
\begin{equation}
\sigma_k \in HF^0(\tau_{V_{k+1}}\cdots\tau_{V_{r-1}}(V_r),\tau_{V_k}\tau_{V_{k+1}}\cdots \tau_{V_{r-1}}(V_r))
\end{equation}
for all $1 \leq k \leq r-1$. Their composition is an element
\begin{equation} \label{eq:composite-t}
\sigma = \sigma_1 \cdots \sigma_{r-1} \in HF^0(V_r,\tau_{V_1}\cdots \tau_{V_{r-1}}(V_r)).
\end{equation}
This all takes place in the fibre $M$. In the total space $E$, consider $L = \Delta_{\gamma_r}$ and its wrapped Floer cohomology $HW^*(L,L)$. The connection between the two is as follows:

\begin{proposition} \label{th:geometric-vanishing}
If $HW^*(L,L)$ is nontrivial, $\sigma$ vanishes.
\end{proposition}
\begin{figure}
\begin{centering}
\begin{picture}(0,0)%
\includegraphics{wrap1.pstex}%
\end{picture}%
\setlength{\unitlength}{3947sp}%
\begingroup\makeatletter\ifx\SetFigFont\undefined%
\gdef\SetFigFont#1#2#3#4#5{%
  \reset@font\fontsize{#1}{#2pt}%
  \fontfamily{#3}\fontseries{#4}\fontshape{#5}%
  \selectfont}%
\fi\endgroup%
\begin{picture}(4645,2920)(31,-2144)
\put(1451,89){\makebox(0,0)[lb]{\smash{{\SetFigFont{11}{13.2}{\rmdefault}{\mddefault}{\updefault}{$b$}%
}}}}
\put(2656,-286){\makebox(0,0)[lb]{\smash{{\SetFigFont{11}{13.2}{\rmdefault}{\mddefault}{\updefault}{$\beta_1$}%
}}}}
\put(1426,539){\makebox(0,0)[lb]{\smash{{\SetFigFont{11}{13.2}{\rmdefault}{\mddefault}{\updefault}{$\beta_0$}%
}}}}
\put(3676,-736){\makebox(0,0)[lb]{\smash{{\SetFigFont{11}{13.2}{\rmdefault}{\mddefault}{\updefault}{$z$}%
}}}}
\end{picture}%
\caption{\label{fig:wrap1}}
\end{centering}
\end{figure}

Before we get to the proof of this, some preliminary considerations are necessary. Let $(L_0,L_1)$ be the Lefschetz thimbles associated to the paths $(\beta_0,\beta_1)$ drawn in Figure \ref{fig:wrap1} (where we assumed that the basis of vanishing paths was the one from Figure \ref{fig:basis2}). Note that $L_1 = L$, while $L_0$ is a compactly supported perturbation of $L^\alpha$ for some angle $\alpha \in (2\pi,4\pi)$. The vanishing paths intersect in two points $(b,z)$. Using the Picard-Lefschetz formula again, one sees that their contributions are
\begin{equation}
\begin{aligned}
& H_b = \Z/2 \text{ placed in degree $0$}, \\
& H_z = HF^{*+n+1}(\tau_{V_1}\cdots \tau_{V_{r-1}}(V_r),V_r).
\end{aligned}
\end{equation}
For the obvious choice of partition $I_0 = \{z\}$, $I_1 = \{b\}$, the only nontrivial differential in the spectral sequence is given by a canonical element of
\begin{equation} \label{eq:duality}
\Hom^*(H_z,H_b[1]) \iso H_{z}^\vee[-1] \iso HF^*(V_r,\tau_{V_1} \cdots \tau_{V_{r-1}}(V_r))
\end{equation}
of degree zero (by Poincar\'e duality in Floer cohomology; recall that here we are dealing with compact Lagrangian submanifolds in the fibres).
\begin{figure}
\begin{centering}
\begin{picture}(0,0)%
\includegraphics{wrap2.pstex}%
\end{picture}%
\setlength{\unitlength}{3947sp}%
\begingroup\makeatletter\ifx\SetFigFont\undefined%
\gdef\SetFigFont#1#2#3#4#5{%
  \reset@font\fontsize{#1}{#2pt}%
  \fontfamily{#3}\fontseries{#4}\fontshape{#5}%
  \selectfont}%
\fi\endgroup%
\begin{picture}(4591,2849)(85,-2073)
\put(1441,164){\makebox(0,0)[lb]{\smash{{\SetFigFont{11}{13.2}{\rmdefault}{\mddefault}{\updefault}{$\tilde{b}$}%
}}}}
\put(3601,-736){\makebox(0,0)[lb]{\smash{{\SetFigFont{11}{13.2}{\rmdefault}{\mddefault}{\updefault}{$\tilde{z}_{2}$}%
}}}}
\put(2666,-286){\makebox(0,0)[lb]{\smash{{\SetFigFont{11}{13.2}{\rmdefault}{\mddefault}{\updefault}{$\tilde{\beta}_1$}%
}}}}
\put(2101,-61){\makebox(0,0)[lb]{\smash{{\SetFigFont{11}{13.2}{\rmdefault}{\mddefault}{\updefault}{$\tilde{z}_1$}%
}}}}
\put(376, 14){\makebox(0,0)[lb]{\smash{{\SetFigFont{11}{13.2}{\rmdefault}{\mddefault}{\updefault}{$\tilde{\beta}_0$}%
}}}}
\end{picture}%
\caption{\label{fig:wrap2}}
\end{centering}
\end{figure}

Instead of computing that element directly, we'll consider the perturbed paths $(\tilde{\beta}_0, \tilde{\beta}_1 = \beta_1)$ from Figure \ref{fig:wrap2}, and their Lefschetz thimbles $(\tilde{L}_0, \tilde{L}_1 = L_1)$. Since this differs from the previous situation by a compactly supported isotopy, we have $HF^*(\tilde{L}_0,\tilde{L}_1) \iso HF^*(L_0,L_1)$. On the other hand, there are now three intersection points $(\tilde{b},\tilde{z}_{1},\tilde{z}_{2})$ in the base, whose contributions are
\begin{equation}
\begin{aligned}
& H_{\tilde{z}_1} = HF^{*+n}(V_r,V_r) \iso H^{*+n}(S^n;\Z/2), \\
& H_{\tilde{b}} = \Z/2 \text{ placed in degree $-n-1$}, \\
& H_{\tilde{z}_2} = H_{z} = HF^{*+n+1}(\tau_{V_1}\cdots \tau_{V_{r-1}}(V_r),V_r). \\
\end{aligned}
\end{equation}
We can choose the filtration to be given by $\tilde{I}_0 = \{\tilde{b},\tilde{z}_{2}\}$ and $\tilde{I}_1 = \{\tilde{z}_1\}$, in which case the only nontrivial differential is a degree $1$ map
\begin{equation} \label{eq:tilde-differential}
H_{\tilde{b}} \oplus H_{\tilde{z}_2} \longrightarrow H_{\tilde{z}_1}.
\end{equation}
The first component is given by an element of $H_{\tilde{z}_1}^{-n} \iso H^0(S^n;\Z/2) = \Z/2$.

\begin{lemma} \label{th:r1}
The first component of \eqref{eq:tilde-differential} is nonzero.
\end{lemma}

\proof Let's return to the proof of Proposition \ref{th:ss}. There is a unique (up to translation) holomorphic strip $w: \R \times [0,1] \rightarrow \C$ which connects $\tilde{b}$ and $\tilde{z}_1$. If we consider the unperturbed almost complex stucture $(J_t)$, the pseudo-holomorphic strips in $E$ contributing to the first component of \eqref{eq:tilde-differential} would have to project to that strip in the base. For the perturbed almost complex structure $(\tilde{J}_t)$ this is still approximately true, meaning that the pseudo-holomorphic strips project to a small neighbourhood of $\overline{w(\R \times [0,1])}$. In particular, for the purposes of computing that component, we can ignore the presence of the other intersection point $\tilde{z}_2$.

However, if there was no other intersection point, we would have an obvious isotopy which merges $\tilde{z}_1$
with $\tilde{b}$. By isotopy invariance, there is necessarily some cancellation between the contributions of these two points, which proves the desired result.
\qed

The second component of \eqref{eq:tilde-differential} is a map $HF^*(\tau_{V_1}\cdots \tau_{V_{r-1}}(V_r),V_r) \rightarrow HF^*(V_r,V_r)$. Following our previous discussion, we prefer to think of it dually as a map
\begin{equation} \label{eq:tilde2}
HF^*(V_r,V_r) \longrightarrow HF^*(V_r,\tau_{V_1}\cdots \tau_{V_{r-1}}(V_r)).
\end{equation}

\begin{lemma}
\eqref{eq:tilde2} is composition with the element $\sigma$ from \eqref{eq:composite-t}.
\end{lemma}

\proof As before, there is a unique holomorphic strip $w$ which connects $\tilde{z}_2$ to $\tilde{z}_1$. Without loss of generality, we may assume that the Lefschetz fibration is symplectically trivial (a product of fibre and base) in a small neighbourhood of each $\tilde{z}_k$. Choose almost complex structures $(J_t)$ which make $\pi_E$ pseudo-holomorphic, and which respect the product structure near the $\tilde{z}_k$. These are of course not generic, but we will see that the moduli spaces appearing in this particular computation can be made regular while remaining within that class; after that, a Gromov compactness argument as in the proof of Proposition \ref{th:ss} applies.

$w$ itself is regular as a holomorphic strip, and its moduli space is zero-dimensional (up to translation). Hence, counting pseudo-holomorphic strips in $E$ which connect the critical points in $E_{\tilde{z}_1}$ to those in $E_{\tilde{z}_2}$ is the same as counting pseudo-holomorphic lifts of $w$ (including regularity; this is by a similar comparison of linearized operators as in the proof of Proposition \ref{th:ss}). Invariants counting pseudo-holomorphic sections of Lefschetz fibrations were defined in \cite{seidel01}. In our case, the fibration is $w^*E \rightarrow \R \times [0,1]$; it can be thought of as being glued together from simpler fibrations over $\R \times [0,1]$, each of which has only one critical point (see Figure \ref{fig:compose} for a schematic picture). The invariants counting pseudo-holomorphic sections of the simpler fibrations are given by composition with the elements $\sigma_k$, essentially by definition of those elements. The TQFT-type structure of the theory (see \cite[Proposition 2.22]{seidel01} for a precise statement of the relevant gluing theorem) then implies that sections of $w^*E \rightarrow [0,1]$ are described by composition with $\sigma$. \qed
\begin{figure}
\begin{centering}
\begin{picture}(0,0)%
\includegraphics{compose.pstex}%
\end{picture}%
\setlength{\unitlength}{3552sp}%
\begingroup\makeatletter\ifx\SetFigFont\undefined%
\gdef\SetFigFont#1#2#3#4#5{%
  \reset@font\fontsize{#1}{#2pt}%
  \fontfamily{#3}\fontseries{#4}\fontshape{#5}%
  \selectfont}%
\fi\endgroup%
\begin{picture}(6277,3024)(211,-2923)
\put(3076,-661){\makebox(0,0)[lb]{\smash{{\SetFigFont{10}{12.0}{\rmdefault}{\mddefault}{\updefault}{$w$}%
}}}}
\put(226,-2011){\makebox(0,0)[lb]{\smash{{\SetFigFont{10}{12.0}{\rmdefault}{\mddefault}{\updefault}{$HF^*(V_{r},V_{r})$}%
}}}}
\put(3901,-2011){\makebox(0,0)[lb]{\smash{{\SetFigFont{10}{12.0}{\rmdefault}{\mddefault}{\updefault}{$HF^*(V_r,\tau_{V_{r-2}}\tau_{V_{r-1}}(V_{r}))$}%
}}}}
\put(1726,-2011){\makebox(0,0)[lb]{\smash{{\SetFigFont{10}{12.0}{\rmdefault}{\mddefault}{\updefault}{$HF^*(V_r,\tau_{V_{r-1}}(V_r))$}%
}}}}
\end{picture}%
\caption{\label{fig:compose}}
\end{centering}
\end{figure}

\proof[Proof of Proposition \ref{th:geometric-vanishing}]
By assumption, $HW^*(L,L)$ is nontrivial. Because of the ring structure, this means that the map $H^*(L;\Z/2) \rightarrow HW^*(L,L)$ which defines the unit in $HW^*(L,L)$ must be nontrivial. One step in this map can be identified with the continuation map
\begin{equation}
\Z/2 \iso HF^*(L_0^{-2\pi},L_1) \longrightarrow HF^*(L_0,L_1).
\end{equation}
From Proposition \ref{th:element}, we know that the image of this map is represented by the unique point of $L_0 \cap L_1$ lying over $b$. This point is always a cocycle, and it is a coboundary if and only if the map \eqref{eq:duality} is nonzero. Hence, by assumption we know that the map is zero, or equivalently, that the total rank of $HF^*(L_0,L_1)$ is one more than that of $HF^{*+n+1}(\tau_{V_1}\cdots \tau_{V_{r-1}}(V_r),V_r)$.

Another way to compute that same rank is to use the perturbed situation from \eqref{eq:tilde-differential}. We know from Lemma \ref{th:r1} that the degree $-n$ element of $H_{\tilde{z}_1}$ always lies in the image of that map. If $\sigma$ is nonzero, the degree $0$ part of \eqref{eq:tilde2} is nonzero, which in view of duality means that the degree $0$ element of $H_{\tilde{z}_1}$ also lies in the image of \eqref{eq:tilde-differential}. Hence, in that case the rank of $HF^*(L_0,L_1)$ is one less than that of $HF^{*+n+1}(\tau_{V_1}\cdots \tau_{V_{r-1}}(V_r),V_r)$, which is a contradiction. The conclusion is that $\sigma$ is zero, as desired. \qed

\section{Algebraic implications\label{sec:algebra}}

We continue to consider the situation from Proposition \ref{th:geometric-vanishing}, and concentrate on the implications of the vanishing of $\sigma$ for the structure of the vanishing cycles $(V_1,\dots,V_r)$. This discussion will take place entirely inside the fibre.

Take the Fukaya category of $M$, whose objects are closed exact Lagrangian submanifolds with gradings. We denote it by $\B$, and also consider its derived category $D(\B) = H^0(\tw(\B))$, defined via twisted complexes as in \cite{kontsevich94,seidel04}. In the derived category, for any two objects $X_0$, $X_1$ one can introduce the twisted object $T_{X_0}(X_1)$, unique up to canonical isomorphism (this is familiar from the theory of mutations \cite{rudakov90}, but the specific notation here is borrowed from \cite[Section 5]{seidel04}). By construction, it comes with a canonical morphism $s_{X_0,X_1} \in \Hom_{D(\B)}(X_1,T_{X_0}(X_1))$. Specialize to $X_0 = V_k$ and $X_1 = T_{V_{k+1}}\cdots T_{V_{r-1}}(V_r)$, and call the resulting morphisms $s_k$. Their composition is an element
\begin{equation}
s = s_1 \cdots s_{r-1} \in \Hom_{D(\B)}(V_r,T_{V_1}\cdots T_{V_{r-1}}(V_{r})).
\end{equation}

\begin{lemma} \label{th:tau-t}
There is an isomorphism $\tau_{V_1}\cdots \tau_{V_{r-1}}(V_{r}) \iso T_{V_1}\cdots T_{V_{r-1}}(V_{r})$ in $D(\B)$, which moreover takes $\sigma$, from \eqref{eq:composite-t}, to $s$.
\end{lemma}

\proof This is a mild extension of material from \cite[Section 17]{seidel04}, and follows from it by diagram-chasing in the following overall structure:
\begin{equation} \label{eq:big-diagram}
\xymatrix{
V_{r} \ar[dr] \ar[r] & \tau_{V_{r-1}}(V_{r}) \ar[r] \ar[dr] & \tau_{V_{r-2}}\tau_{V_{r-1}}(V_{r}) \cdots \\
& T_{V_{r-1}}(V_{r}) \ar[u] \ar[dr] & T_{V_{r-2}}(\tau_{V_{r-1}}(V_{r})) \ar[u] \cdots \\
& & T_{V_{r-2}}(T_{V_{r-1}}(V_{r})) \cdots \ar[u]
}
\end{equation}
Consider first the leftmost triangle
\begin{equation}
\xymatrix{
V_{r} \ar[dr]_-{s_{r-1}} \ar[r]^-{\sigma_{r-1}} & \tau_{V_{r-1}}(V_r) \\ & T_{V_{r-1}}(V_r) \ar[u]_-{\iso}
}
\end{equation}
The $\uparrow$ is the isomorphism from \cite[Theorem 17.16]{seidel04}, which by construction makes the triangle commute. The same holds for all other triangles in \eqref{eq:big-diagram}. Now consider one of the lozenges,
\begin{equation} \label{eq:lozenge}
\xymatrix{
\tau_{V_{r-1}}(V_{r}) \ar[dr] & \\ T_{V_{r-1}}(V_{r}) \ar[u]^-{\iso} \ar[dr]_{s_{r-2}} & T_{V_{r-2}}(\tau_{V_{r-1}}(V_{r})) \\ & T_{V_{r-2}}(T_{V_{r-1}}(V_{r})) \ar[u]_-{T_{V_{r-2}}(\iso)}
}
\end{equation}
Here, the $\searrow$'s are the canonical elements $s_{X_0,X_1} \in \Hom_{D(\B)}(X_1,T_{X_0}(X_1))$ for $(X_0,X_1) = (V_{r-2},\tau_{V_{r-1}}(V_{r}))$ and $(X_0,X_1) = (V_{r-2},T_{V_{r-1}}(V_{r}))$ respectively, which are purely algebraic. It is clear from the definition that $T_{V_{r-2}}$ is actually an exact functor from $D(\B)$ to itself, and that the canonical elements form a natural transformation from the identity functor to $T_{V_{r-2}}$. Now, the $\uparrow$'s in \eqref{eq:lozenge} are, respectively, the isomorphism from \cite[Theorem 17.16]{seidel04} and its image under $T_{V_{r-2}}$. Hence, the diagram commutes by naturality, and the same argument applies to the rest of \eqref{eq:big-diagram}. By going around the sides of that, one obtains the desired result. \qed

Let $\A$ be the directed $A_\infty$-subcategory \cite[Section 5n]{seidel04} associated to the objects $(V_1,\dots,V_{r-1})$ in $\B$. This comes with a canonical (up to quasi-isomorphism) $A_\infty$-functor $\A \rightarrow \B$, which then induces an exact functor $D(\A) \rightarrow D(\B)$.

\begin{proposition} \label{th:summand}
Suppose that $s = 0$. Then $V_r$ is isomorphic to a direct summand of an object lying in the image of the functor $D(\A) \rightarrow D(\B)$.
\end{proposition}

\proof Choose a cochain representative of $s$ in $\hom_{\tw(\B)}^0(V_{r},T_{V_1}\cdots T_{V_{r-1}}(V_{r}))$, and let $C$ be its mapping cone. By construction, this fits into a distinguished triangle in $D(\B)$ of the form
\begin{equation}
\xymatrix{V_r \ar[r]^-{s} & T_{V_1}\cdots T_{V_{r-1}}(V_{r}) \ar[d] \\ & C \ar[ul]^-{[1]}}
\end{equation}
Our assumption implies that the triangle splits, which means that
\begin{equation} \label{eq:summand}
V_{r}[1] \oplus T_{V_1}\cdots T_{V_{r-1}}(V_{r}) \iso C.
\end{equation}
It remains to write down $C$ more explicitly. For this purely algebraic question, we find it convenient to replace $\B$ by a quasi-isomorphic $A_\infty$-category $\tilde\B$, with the same objects, which is strictly unital. That can always be done, see \cite[Section 2]{seidel04} or \cite[Theorem 3.2.1.1]{lefevre}. If $\tilde\A$ is the directed $A_\infty$-subcategory associated to $(V_1,\dots,V_{r-1})$ in $\tilde\B$, we have a diagram (commutative up to isomorphism)
\begin{equation}
\xymatrix{D(\A) \ar[d]^-{\iso} \ar[r] & D(\B) \ar[d]^-{\iso} \\
D(\tilde\A) \ar[r] & D(\tilde\B)}
\end{equation}
where the vertical arrows are equivalences. Hence, the essential situation does not change, but we do get two technical simplifications. First, the $A_\infty$-functor $\tilde\A \rightarrow \tilde\B$ becomes an embedding of an $A_\infty$-subcategory. Second, there are simple canonical representatives for the objects $T_{X_0}(X_1)$ in $\tw(\tilde\B)$, namely
\begin{equation} \label{eq:t-cone}
T_{X_0}(X_1) = \mathit{Cone}(\hom_{\tw(\tilde\B)}(X_0,X_1) \otimes X_0 \rightarrow X_1),
\end{equation}
and the canonical morphism $s_{X_0,X_1}$ is then represented by the obvious inclusion of twisted complexes, $X_1 \hookrightarrow T_{X_0}(X_1)$. Iterating this picture, one finds that $s$ is represented by the inclusion $V_r \hookrightarrow T_{V_1}\cdots T_{V_{r-1}}(V_r)$. The mapping cone is then the quotient $\tilde{C}$ of that inclusion, which is of the form
\begin{equation} \label{eq:tilde-c}
\tilde{C} = \bigoplus \hom_{\tilde\B}(V_{k_{i-1}},V_{k_i})[1] \otimes \cdots \otimes \hom_{\tilde\B}(V_{k_0},V_{k_1})[1] \otimes V_{k_0},
\end{equation}
where the sum is over all $i \geq 1$ and $1 \leq k_0 < k_1 < \cdots < k_i = r$. The crucial fact about $\tilde{C}$, which is a direct consequence of its definition \eqref{eq:t-cone}, is the following directedness property:
\begin{equation}
\parbox{30em}{
The differential $\partial_{\tilde{C}}$ has no nontrivial entries which decrease $k_0$. Moreover, the entries which preserve $k_0$ are of the form $\phi_{k_0} \otimes e_{V_{k_0}}$, where $\phi_{k_0}$ is an endomorphism of the vector space $\hom_{\tilde\B}(V_{k_{i-1}},V_{k_i})[1] \otimes \cdots \otimes \hom_{\tilde\B}(V_{k_0},V_{k_1})[1]$, and $e_{V_{k_0}}$ is the strict identity morphism.
}
\end{equation}
Hence, $\tilde{C}$ is in fact an object of $\tw(\tilde{\A}) \subset \tw(\tilde{\B})$. On the other hand, by construction it corresponds to $C$ under the equivalence of triangulated categories $D(\tilde\B) \iso D(\B)$. A look at \eqref{eq:summand} (or rather its version shifted by $-1$) then yields the desired result. \qed

\begin{corollary} \label{th:algebraic-vanishing}
If $HW^*(\Delta_{\gamma_r},\Delta_{\gamma_r}) \neq 0$, then $V_r$ is isomorphic to a direct summand of an object lying in the image of the functor $D(\A) \rightarrow D(\B)$.
\end{corollary}

This follows directly from Proposition \ref{th:geometric-vanishing}, Lemma \ref{th:tau-t}, and Proposition \ref{th:summand}.

\section{$(A_m)$ Milnor fibres\label{sec:milnor}}

Fix $m \geq 1$, $n \geq 2$. The $2n$-dimensional $(A_m)$ type Milnor fibre is the affine hypersurface
\begin{equation}
M_m = \{ x_1^2 + \cdots + x_n^2 + x_{n+1}^{m+1} = 1\} \subset \C^{n+1}.
\end{equation}
When equipped with the restriction of the standard symplectic form, this is a Liouville manifold (and like any affine hypersurface, comes with a canonical trivialization of its anticanonical bundle). Let $\sqrt[m+1]{1} \subset \C$ be the subset of $(m+1)$-st roots of unity. To every embedded path $\delta$ whose endpoints lie in $\sqrt[m+1]{1}$, and which avoids that set otherwise, one can associate a Lagrangian sphere $S_\delta \subset M_m$ \cite[Section 6c]{khovanov-seidel98}.

\begin{lemma} \label{th:twist-twist}
Let $\delta_0,\delta_1$ be any two paths. Denote by $t_{\delta_0}$ the right handed half-twist, which is a diffeomorphism of $\C$ preserving $\sqrt[m+1]{1}$, and by $\tau_{S_{\delta_0}}$ the Dehn twist along the associated Lagrangian sphere. Then there is an isotopy of Lagrangian spheres,
\begin{equation} \label{eq:twist-twist}
S_{t_{\delta_0}(\delta_1)} \htp \tau_{S_{\delta_0}}(S_{\delta_1}).
\end{equation}
\end{lemma}

\proof[Sketch of proof] There is a canonical symplectic fibration over configuration space $\mathit{Conf}_{m+1}(\C)$ with fibre $M_m$. The Picard-Lefschetz theorem shows that $\tau_{S_{\delta_0}}$ is the monodromy of that fibration along a particular loop in the base. That path corresponds to a braid, and if one realizes that braid as diffeomorphism of $(\C,\sqrt[m+1]{1})$, the result is precisely $t_{\delta_0}$. By combining these two facts, one can explicitly construct a family of Lagrangian spheres interpolating between the two sides of \eqref{eq:twist-twist}. Compare the discussion in \cite[Remark 16.14]{seidel04}. \qed

\begin{lemma} \label{th:non-isotopic}
Let $\delta_0,\delta_1$ be two paths which are not isotopic (within the class of paths we've been considering). Then the image of the product map
\begin{equation} \label{eq:triangle}
HF^*(S_{\delta_1},S_{\delta_0}) \otimes HF^*(S_{\delta_0},S_{\delta_1}) \longrightarrow HF^*(S_{\delta_0}, S_{\delta_0}) \iso H^*(S_{\delta_0};\Z/2).
\end{equation}
is contained in $H^n(S_{\delta_0};\Z/2)$.
\end{lemma}
\begin{figure}
\begin{centering}
\begin{picture}(0,0)%
\includegraphics{intersect.pstex}%
\end{picture}%
\setlength{\unitlength}{3355sp}%
\begingroup\makeatletter\ifx\SetFigFont\undefined%
\gdef\SetFigFont#1#2#3#4#5{%
  \reset@font\fontsize{#1}{#2pt}%
  \fontfamily{#3}\fontseries{#4}\fontshape{#5}%
  \selectfont}%
\fi\endgroup%
\begin{picture}(3499,3169)(1239,-1868)
\put(4501, 14){\makebox(0,0)[lb]{\smash{{\SetFigFont{10}{12.0}{\rmdefault}{\mddefault}{\updefault}{$\xi$}%
}}}}
\put(3751,-1186){\makebox(0,0)[lb]{\smash{{\SetFigFont{10}{12.0}{\rmdefault}{\mddefault}{\updefault}{$\delta_0$}%
}}}}
\put(3076,164){\makebox(0,0)[lb]{\smash{{\SetFigFont{10}{12.0}{\rmdefault}{\mddefault}{\updefault}{$\delta_1$}%
}}}}
\end{picture}%
\caption{\label{fig:intersect}}
\end{centering}
\end{figure}

\proof[Sketch of proof] Without loss of generality, suppose that $\delta_1$ is as Figure \ref{fig:intersect}. Take the infinite path $\xi$ and associate to it a properly embedded Lagrangian submanifold $L \subset M$, $L \iso \R \times S^{n-1}$, following the same construction as for the spheres $S_\delta$. Since $\delta_0$ is not isotopic to $\delta_1$, it has essential intersection with $\xi$, so their geometric intersection number is $I(\delta_0,\xi) > 0$. By the same argument as in \cite[Lemma 6.19]{khovanov-seidel98} we have
\begin{equation} \label{eq:intersection-number}
\begin{aligned}
& \mathrm{dim}\, HF^*(L,S_{\delta_0}) = 2 I(\delta_0,\xi), \\
& HF^*(L,S_{\delta_1}) = 0.
\end{aligned}
\end{equation}
Suppose that there are elements $a_2 \in HF^*(S_{\delta_1},S_{\delta_0})$ and $a_1 \in HF^*(S_{\delta_0},S_{\delta_1})$ whose product does not lie in $H^n(S_{\delta_0};\Z/2) \subset HF^*(S_{\delta_0},S_{\delta_0})$. This product is then necessarily an invertible element of the ring $HF^*(S_{\delta_0},S_{\delta_0}) \iso H^*(S^n;\Z/2)$, which means that
\begin{equation}
HF^*(L,S_{\delta_0}) \xrightarrow{a_1 \cdot} HF^*(L,S_{\delta_1}) \xrightarrow{a_2 \cdot} HF^*(L,S_{\delta_0})
\end{equation}
is an isomorphism, contradicting \eqref{eq:intersection-number}. \qed
\begin{figure}
\begin{centering}
\begin{picture}(0,0)%
\includegraphics{clockwise.pstex}%
\end{picture}%
\setlength{\unitlength}{3355sp}%
\begingroup\makeatletter\ifx\SetFigFont\undefined%
\gdef\SetFigFont#1#2#3#4#5{%
  \reset@font\fontsize{#1}{#2pt}%
  \fontfamily{#3}\fontseries{#4}\fontshape{#5}%
  \selectfont}%
\fi\endgroup%
\begin{picture}(2313,2370)(1336,-1859)
\put(2476,-1786){\makebox(0,0)[lb]{\smash{{\SetFigFont{10}{12.0}{\rmdefault}{\mddefault}{\updefault}{$\delta_4$}%
}}}}
\put(3001,239){\makebox(0,0)[lb]{\smash{{\SetFigFont{10}{12.0}{\rmdefault}{\mddefault}{\updefault}{$\delta_3$}%
}}}}
\put(1351,-961){\makebox(0,0)[lb]{\smash{{\SetFigFont{10}{12.0}{\rmdefault}{\mddefault}{\updefault}{$\delta_1$}%
}}}}
\put(2851,-1036){\makebox(0,0)[lb]{\smash{{\SetFigFont{10}{12.0}{\rmdefault}{\mddefault}{\updefault}{$\delta_2$}%
}}}}
\end{picture}%
\caption{\label{fig:clockwise}}
\end{centering}
\end{figure}

Write $\delta^{k,l}$ for the straight line segment connecting $e^{2\pi i k/(m+1)}$ to $e^{2 \pi i l/(m+1)}$, where $k \neq l$ mod $m+1$. Suppose first that we choose vanishing cycles $V_j = S_{\delta^{j-1,j}}$ for $j = 1,\dots,m$. Then the Liouville $(2n+2)$-manifold $E$ constructed from $M_m$ and $(V_1,\dots,V_m)$ is Liouville isomorphic to standard symplectic $\R^{2n+2}$. This is just the Morsification of the singularity of type $(A_m)$, explained in reverse. In general, we can apply Hurwitz moves to a given collection, and this gives new collections of vanishing cycles which still lead to standard $\R^{2n+2}$ as the total space. In our particular case, the collections $(\tilde{V}_1,\dots,\tilde{V}_m)$ obtained by applying Hurwitz moves to $(V_1,\dots,V_m)$ are precisely those of the following kind. Each $\tilde{V}_j = S_{\tilde{\delta}_j}$ is associated to some straight line segment $\tilde\delta_j = \delta^{k_j,l_j}$. Moreover, any two such segments are either disjoint or intersect at a single endpoint, and the union of all of them forms a tree inside the unit disc. Finally, if several segments meet at a common endpoint, the directions at that point are in clockwise order (see Figure \ref{fig:clockwise} for an example).
An obvious consequence is that only finitely many different collections $(\tilde{V}_1,\dots,\tilde{V}_m)$ arise (exactly $(m+1)^{m-1}$, which is Cayley's formula for the number of trees with $m$ numbered edges and with an additional choice of distinguished vertex).

\begin{remark}
There is a more geometric way of seeing how that particular number arises (this is not new, compare for instance  \cite[Introduction]{arnold96}). Let ${\mathcal P}_m$ be the space of all polynomials of degree $m+1$ which are monic, have zero subleading term, and moreover have $m$ distinct critical values. The Lyashko-Looijenga map
\begin{equation}
{\mathcal P}_m \longrightarrow \mathit{Conf}_m(\C)
\end{equation}
which associates to each such polynomial its critical values, is a covering of degree $(m+1)^{m-1}$ (see \cite{looijenga74}, or \cite[Chapter 5]{lando-zvonkin} for an expository account). We'd like to view ${\mathcal P}_m$ as the space of all Lefschetz fibrations $\C^{n+1} \rightarrow \C$ of the form $x \mapsto x_1^2 + \cdots + x_n^2 + p(x_{n+1})$ (here, Lefschetz fibration is understood in a slightly looser sense than in Section \ref{sec:lefschetz}, so as to fit into the natural algebro-geometric framework). Fix a base point in ${\mathcal P}_m$, and a basis of vanishing paths for the Lefschetz fibration associated to that point, and consider the resulting collection of vanishing cycles. Going around any loop in $\mathit{Conf}_m(\C)$ transforms this into another basis of vanishing paths, to which corresponds an a priori different collection of vanishing cycles (this is a version of the braid group action by Hurwitz moves). However, if the loop can be lifted to ${\mathcal P}_m$, then the new collection of vanishing cycles is isotopic to the previous one, by a monodromy argument. This immediately shows that there are most $(m+1)^{m-1}$ different such collections (however, additional work is required to show that there are not less than that).
\end{remark}

\begin{lemma} \label{th:standard}
Choose $V_j = S_{\delta^{j-1,j}}$ for $j = 1,\dots,m$, and $V_{m+1} = S_{\delta^{k,l}}$ for any $k \neq l$ mod $m+1$. Then the Liouville manifold $E$ constructed from $M_m$ and $(V_1,\dots,V_{m+1})$ is Liouville isomorphic to standard symplectic $T^*\!S^{n+1}$.
\end{lemma}

\proof[Sketch of proof] Let's consider first the toy model case $m = 1$, $V_1 = V_2 = S_{\delta^{0,1}}$. What we have is a fibre $M_1 \iso T^*\!S^n$ and two vanishing cycles which are both equal to the zero-section $S^n$, which is the standard Lefschetz fibration with total space $T^*\!S^{n+1}$; see for instance \cite[Section 5]{maydanskiy09}. A similar argument works for all $m$ if $(k,l) = (m-1,m)$. In that situation, one has two equal vanishing cycles forming a Lagrangian sphere in $E$, while the remaining vanishing cycles provide handle attachments that cancel out the extra topology of the fibre.
\begin{figure}
\begin{centering}
\begin{picture}(0,0)%
\includegraphics{snake.pstex}%
\end{picture}%
\setlength{\unitlength}{3355sp}%
\begingroup\makeatletter\ifx\SetFigFont\undefined%
\gdef\SetFigFont#1#2#3#4#5{%
  \reset@font\fontsize{#1}{#2pt}%
  \fontfamily{#3}\fontseries{#4}\fontshape{#5}%
  \selectfont}%
\fi\endgroup%
\begin{picture}(2538,2970)(4411,-2234)
\put(4951,539){\makebox(0,0)[lb]{\smash{{\SetFigFont{10}{12.0}{\rmdefault}{\mddefault}{\updefault}{$\tilde{\delta}_3$}%
}}}}
\put(6226,-811){\makebox(0,0)[lb]{\smash{{\SetFigFont{10}{12.0}{\rmdefault}{\mddefault}{\updefault}{$\tilde{\delta}_7$}%
}}}}
\put(5626,-436){\makebox(0,0)[lb]{\smash{{\SetFigFont{10}{12.0}{\rmdefault}{\mddefault}{\updefault}{$\tilde{\delta}_8$}%
}}}}
\put(5101,-1711){\makebox(0,0)[lb]{\smash{{\SetFigFont{10}{12.0}{\rmdefault}{\mddefault}{\updefault}{$\tilde{\delta}_6$}%
}}}}
\put(4426,-1036){\makebox(0,0)[lb]{\smash{{\SetFigFont{10}{12.0}{\rmdefault}{\mddefault}{\updefault}{$\tilde{\delta}_5$}%
}}}}
\put(4426,-61){\makebox(0,0)[lb]{\smash{{\SetFigFont{10}{12.0}{\rmdefault}{\mddefault}{\updefault}{$\tilde{\delta}_4$}%
}}}}
\put(4651,-2161){\makebox(0,0)[lb]{\smash{{\SetFigFont{10}{12.0}{\rmdefault}{\mddefault}{\updefault}{$m = 7$ and $(k,l) = (1,6)$}%
}}}}
\put(6751,-1111){\makebox(0,0)[lb]{\smash{{\SetFigFont{10}{12.0}{\rmdefault}{\mddefault}{\updefault}{$\tilde{\delta}_1$}%
}}}}
\put(6751,-61){\makebox(0,0)[lb]{\smash{{\SetFigFont{10}{12.0}{\rmdefault}{\mddefault}{\updefault}{$\tilde{\delta}_2$}%
}}}}
\end{picture}%
\caption{\label{fig:snake}}
\end{centering}
\end{figure}

In the general case one can use Hurwitz moves, applied only to the first $m$ cycles, to modify the given collection to $\tilde{V}_j = S_{\tilde{\delta}_j}$, where
\begin{equation} \label{eq:new-collection}
\begin{aligned}
& \tilde{\delta}_1 = \delta^{l+1,l+2},\dots, \tilde{\delta}_{m+k-l} = \delta^{m+k,m+1+k}, \\
& \tilde{\delta}_{m+k-l+1} = \delta^{k+1,k+2},\dots,\tilde{\delta}_{m-1} = \delta^{l-1,l}, \\
& \tilde{\delta}_m = \tilde{\delta}_{m+1} = \delta^{k,l}.
\end{aligned}
\end{equation}
Here, we have assumed without loss of generality that $0 \leq k < l \leq m$ (see Figure \ref{fig:snake} for a picture of the paths $\tilde\delta_j$). Explicitly, if $\sigma_i$ denotes the $i$-th elementary Hurwitz move in the conventions from \cite[Section 16]{seidel04}, then the modification leading to \eqref{eq:new-collection} is given by $\sigma_{m-1}^{-1}\sigma_{m-2}^{-1} \cdots \sigma_{k-l+m+1}^{-1} (\sigma_1\sigma_2 \cdots \sigma_{m-1})^{m-l}$, where the order is from right to left (alternatively, the existence of such a sequence of moves can be derived from the general discussion preceding this Lemma). This gets us back to a version of the previous situation, since the two last vanishing cycles coincide, while the others again contribute cancelling handle attachments. \qed

\section{$(A_m)$ quiver representations\label{sec:quiver}}

Consider the directed quiver of type $(A_m)$, for some $m \geq 1$. A representation of this quiver is a sequence of finite-dimensional vector spaces $W_i$, $1 \leq i \leq m$, and linear maps $\rho_i: W_i \rightarrow W_{i+1}$, $1 \leq i \leq m-1$. This can be over an arbitrary field, but the relevant case for us is where the ground field is $\Z/2$. An elementary case of Gabriel's theorem \cite{gabriel62} says that any indecomposable representation of the $(A_m)$ quiver is isomorphic to one of the following form:
\begin{equation}
W_i = \begin{cases}
\Z/2 & k < i \leq l, \\
0 & \text{otherwise.}
\end{cases}, \qquad
\rho_i = \begin{cases}
1 & k < i < l, \\
0 & \text{otherwise.}
\end{cases}
\end{equation}
Here, $0 \leq k < l \leq m$, so there are a total of $m(m+1)/2$ different indecomposable representations.

We will now recast this statement as one about twisted complexes over certain $A_\infty$-categories. Consider the $A_\infty$-category $\A_m$ over $\Z/2$, which has $m$ objects denoted by $(V_1,\dots,V_m)$, is strictly unital, and has morphism spaces
\begin{equation} \label{eq:am-quiver}
\mathit{hom}_{\A_m}(V_i,V_j) = \begin{cases} \Z/2 \cdot e_i & \text{for $i = j$, where $e_i$ is the unit}, \\
\Z/2 \cdot f_i & \text{for $i = j-1$, where $f_i$ has degree $1$}, \\
0 & \text{otherwise.}
\end{cases}
\end{equation}
This determines the $A_\infty$-structure of $\A_m$ uniquely: the only nonvanishing products are $\mu^2(e_i,e_i) = e_i$ as well as $\mu^2(f_i,e_i) = f_i = \mu^2(e_{i+1},f_i)$. Objects of $\tw(\A_m)$ are generally formal sums
\begin{equation} \label{eq:directed-twisted}
C = \bigoplus_{i=1}^m W_i \otimes V_i
\end{equation}
where the $W_i$ are finite-dimensional graded vector spaces over $\Z/2$, together with a differential which is a formal matrix $\partial_C = (\partial_{C,ji})$ consisting of
\begin{equation}
\partial_{C,ji} \in \big(\hom_{\Z/2}(W_i,W_j) \otimes \hom_{\A_m}(V_i,V_j)\big)^1.
\end{equation}
Without changing the quasi-equivalence type of $\tw(\A_m)$, one can restrict to twisted complexes where $\partial_{C,ii} = 0$. This is a general fact about directed $A_\infty$-categories \cite[Remark 5.19]{seidel04}. For the specific case of $\A_m$, this means that the only nonzero components of $\partial$ are $\partial_{C,i+1,i} = \rho_i \otimes f_i$, where $\rho_i \in \hom_{\Z/2}(W_i, W_{i+1})$ is a degree $0$ linear map. Hence $C$ splits as a direct sum corresponding to the graded pieces of the associated vector spaces. Moreover, each such piece is precisely given by a representation of the $(A_m)$ quiver. In particular, we can consider the twisted complexes $C^{k,l}$ corresponding to the indecomposable representations considered above, and as an immediate consequence,

\begin{lemma} \label{th:gabriel}
Every indecomposable object of $D(\A_m) = H^0(\tw(\A_m))$ is isomorphic to a shifted version of $C^{k,l}$, for some $0 \leq k<l \leq m$. \qed
\end{lemma}

\begin{remark}
Since $\A_m$ is directed, the objects $V_i$ (or $C^{i-1,i}$, which is the same) form a full exceptional collection in $H^0(\tw(\A_m))$. There is an action of $\mathit{Br}_m$ on isomorphism classes of such collections, by mutations \cite{rudakov90}. Any object that appears in a mutated collection $(\tilde{V}_1,\dots,\tilde{V}_m)$ is indecomposable, and therefore isomorphic to a shifted version of some $C^{k,l}$. In particular, up to isomorphisms and shifts, only finitely many different exceptional collections arise through the mutation process. To relate this to the geometric finiteness phenomena from the previous section, we note that $\A_m$ is the directed Fukaya category associated to the Lefschetz fibration with fibre $M_m$ and total space $\R^{2n+2}$. By \cite[Corollary 17.17]{seidel04}, mutation of exceptional collections in $H^0(\tw(\A_m))$ corresponds to Hurwitz moves on vanishing cycles. This shows that exactly $(m+1)^{m-1}$ essentially different exceptional collections arise.
\end{remark}

\section{The construction\label{sec:main}}

As outlined in the Introduction, fix some $m,n \geq 2$ and the corresponding Milnor fibre $M_m$. Take vanishing cycles $V_j = S_{\delta_j}$, where $\delta_j = \delta^{j-1,j}$ for $j \leq m$, while $\delta_{m+1}$ can be arbitrary. Form the associated Liouville $(2n+2)$-manifold $E$.

\proof[Proof of Lemma \ref{th:isotopy}] It was shown in \cite[Section 5]{maydanskiy09} that for even $n$, the isotopy class of $S_\delta \subset M_m$ depends only on the endpoints of $\delta$. In fact, the isotopies constructed there are through totally real submanifolds. The diffeomorphism class of $E$ and the homotopy class of its almost complex structure are preserved under such an isotopy. Since any $\delta_{m+1}$ has the same endpoints as some $\delta^{k,l}$, the result follows from Lemma \ref{th:standard}. \qed

\proof[Proof of Lemma \ref{th:haefliger}] We can first attach handles corresponding to $(V_1,\dots,V_m)$, which produces a standard ball $B^{2n+2}$, and then attach a final handle, whose attaching sphere is a Legendrian embedding $S^n \hookrightarrow S^{2n+1} = \partial B^{2n+2}$ derived from $V_{m+1}$. From the topological viewpoint, the data that matter are the isotopy class of the sphere and its framing (trivialization of the normal bundle).

For any $n \geq 3$, any two embeddings $S^n \hookrightarrow S^{2n+1}$ are differentiably isotopic \cite{haefliger62}. Hence, we may assume that the attaching sphere is standard. In that case, an equivalent picture is that $E$ is obtained from two copies of $B^{n+1} \times \R^{n+1}$ by identifying the boundaries through a fibrewise linear automorphism of $S^n \times \R^{n+1}$. In other words, $E$ is the total space of a rank $(n+1)$ vector bundle $\eta \rightarrow S^{n+1}$. Such bundles are classified by the homotopy class of their clutching functions, lying in $\pi_n(O(n+1))$, which is equivalent to the framing data in our previous picture. For $n$ odd, that group sits in a short exact sequence
\begin{equation} \label{eq:split}
0 \rightarrow \pi_{n+1}(S^{n+1}) \rightarrow \pi_n(O(n+1)) \rightarrow \pi_n(O(\infty)) \rightarrow 0.
\end{equation}
The image of a class in $\pi_n(O(\infty))$ determines the stable isomorphism type of $\eta$. Since $TS^{n+1}$ is stably trivial, this is the same as the stable isomorphism type of $TE|S^{n+1}$. On the other hand, we have the natural map $\pi_n(O(n+1)) \rightarrow \pi_n(S^n)$, which computes the Euler class of $\eta$, or equivalently the selfintersection number of $S^{n+1} \subset E$. Since the composition $\pi_{n+1}(S^{n+1}) \rightarrow \pi_n(O(n+1)) \rightarrow \pi_n(S^n)$ is multiplication by $2$ (this is part of the standard computation of the first nontrivial homotopy groups of Stiefel manifolds, see \cite[\S 25]{steenrod} or \cite[Section 8.11]{husemoller}), the selfintersection number detects the left hand subgroup in \eqref{eq:split}.

By definition, $M = M_m$ is a smooth affine hypersurface in $\C^{n+1}$. Hence, its tangent bundle is stably trivial, $TM \oplus \C \iso M \times \C^{n+1}$. Moreover, each of our spheres $V_i$ bounds a Lagrangian ball in $\C^{n+1}$ (because they can be constructed as vanishing cycles for a Lefschetz fibration $\C^{n+1} \rightarrow \C$ with fibre $M$). Hence, we also have a stable trivialization $TV_i \oplus \R \iso V_i \times \R^{n+1}$, and this is compatible with the canonical isomorphism $TM|V_i \iso TV_i \otimes_\R \C$. As a result of this, the manifold $E$ again comes with a stable trivialization $TE \oplus \C \iso E \times \C^{n+2}$. In fact, for dimension reasons this implies that $TE$ itself is trivial. As a consequence of this and the previous discussion, the only topological invariant that can distinguish different $E$'s is the selfintersection number.

Computing the intersection pairing on the total space of a Lefschetz fibration is a standard exercise. Let $H_{n+1}(E,M)$ be the homology of the total space relative to a fibre at infinity. This is generated by the classes of the Lefschetz thimbles $\Delta_1,\dots,\Delta_{m+1}$, and carries a non-symmetric extension of the intersection pairing, which we denote by $\circ$ (in singularity theory, this appears as the linking pairing on the Milnor fibre, see \cite[\S 6]{lamotke75}). In our case, writing $[V_{m+1}] = \sum_{i=1}^m a_i [V_i]$ we have
\begin{equation}
\begin{cases}
\Delta_i \circ \Delta_i = \sigma, & \\
\Delta_i \circ \Delta_j = (-1)^{n+1} V_i \cdot V_j & \text{for $i < j$}, \\
\Delta_i \circ \Delta_j = 0 & \text{for $i > j$.}
\end{cases}
\end{equation}
where $\sigma = (-1)^{\half (n+1)(n+2)}$. The map $H_{n+1}(E) \rightarrow H_{n+1}(E,M)$ takes the generator $x$ to $[\Delta_{m+1}] - \sum_{i=1}^m a_i [\Delta_i]$ and is compatible with the intersection pairing, hence
\begin{equation}
x \cdot x = \sigma (1 - \sum_i a_ia_{i-1} + \sum_i a_i^2).
\end{equation}
On $T^*\!S^{n+1}$ we have $x \cdot x = \sigma \chi(S^{n+1}) = 2\sigma$, hence the condition that our selfintersection number should be the same translates to
\begin{equation}
2 \sum_i a_i^2 - \sum_i a_i a_{i-1} - \sum_i a_i a_{i+1} = 2.
\end{equation}
The left hand side is the standard $(A_m)$ quadratic form, which is positive definite. The only elements which take value $2$ are $a = (0,\dots,0,1,\dots,1,0,\dots,0)$, matching the condition imposed in the statement of the Lemma. \qed

\proof[Proof of Theorem \ref{th:main}] As in Section \ref{sec:algebra}, let $\B$ be the Fukaya category of $M_m$, and $\A$ the directed $A_\infty$-subcategory associated to the collection $(V_1,\dots,V_m)$. The choice of paths means that for $i < m$, $V_i$ intersects $V_{i+1}$ transversally in a single point. We can choose the gradings of the $V_k$ in such a way that the unique generator of $HF^*(V_j,V_{j+1})$ has degree $1$. Since $V_i \cap V_j = \emptyset$ for all $i,j \leq m$ with $|i-j| \geq 2$, a comparison with \eqref{eq:am-quiver} shows that $\A$ is isomorphic to $\A_m$.

By repeatedly applying Lemma \ref{th:twist-twist} and \cite[Theorem 17.16]{seidel04}, one sees that in $D(\B)$,
\begin{equation}
S_{\delta^{k,l}} \iso \tau_{V_{k+1}} \cdots \tau_{V_{l-1}}(V_l) \iso
T_{V_{k+1}} \cdots T_{V_{l-1}}(V_l).
\end{equation}
After writing out explicitly the right hand side (which is easy to do by induction on $l-k$), one sees that it is precisely the image of $C^{k,l}$ under the functor $D(\A) \rightarrow D(\B)$.

Suppose that $E$ contains a Lagrangian sphere representing a nonzero element of $H_{n+1}(E)$. It follows from the handle attachment picture that $\Delta_{m+1}$ is the dual generator of $H^{n+1}(E)$, hence its intersection number with our sphere is nonzero. In view of Lemma \ref{th:intersection-number} and Corollary \ref{th:algebraic-vanishing}, this implies that in $D(\B)$, $V_{m+1}$ is a direct summand of an object $C$ lying in the image of $D(\A) \rightarrow D(\B)$. This means that the product
\begin{multline} \label{eq:product}
\Hom_{D(\B)}(C,V_{m+1}) \otimes \Hom_{D(\B)}(V_{m+1},C) \longrightarrow \Hom_{D(\B)}(V_{m+1},V_{m+1}) \\ = HF^*(V_{m+1},V_{m+1}) \iso H^*(S^n;\Z/2)
\end{multline}
contains the identity in its image. By Lemma \ref{th:gabriel}, $C$ is necessarily a direct sum of shifted copies of various $C^{k,l}$'s, which geometrically means a direct sum of copies of the $S_{\delta^{k,l}}$ with various gradings. But then, our statement concerning \eqref{eq:product} contradicts Lemma \ref{th:non-isotopic}, unless $\delta_{m+1}$ is isotopic to one of the $\delta^{k,l}$.

The conclusion is that if $\delta_{m+1}$ is not isotopic to any $\delta^{k,l}$, then $E$ can't contain a Lagrangian sphere which is nontrivial in homology. Hence, it's not symplectically isomorphic to $T^*\!S^{n+1}$. In the other direction, we already know from Lemma \ref{th:standard} that for any of the $\frac{1}{2} m(m+1)$ choices $\delta_{m+1} = \delta^{k,l}$, the resulting $E$ is isomorphic to $T^*\! S^{n+1}$. \qed


\end{document}


\setcounter{section}{1}
\maketitle

In the proof of \cite[Lemma 1.1]{maydanskiy-seidel09}, we appealed to an explicit isotopy of totally real spheres, constructed in \cite[Section 5]{maydanskiy09b}. That construction works in the lowest dimension ($n = 2$), but is wrong in general (one of the endpoints is not the desired sphere). Here, we explain a different approach, leading to a corrected version of \cite[Lemma 1.1]{maydanskiy-seidel09}, which requires an additional assumption. Independently, while \cite[Lemma 1.2]{maydanskiy-seidel09} makes a statement about homotopy classes of almost complex structures, its proof only determines the isomorphism class of the tangent bundle as an abstract complex vector bundle, which is a priori a weaker statement. The argument here also fills that gap. The rest of the original paper is unaffected.

Consider $M = M_m$ as in \cite{maydanskiy-seidel09}, in complex dimension $n>2$. The construction of $E$ depends on a choice of Lagrangian sphere $S = S_{\delta_{m+1}} \subset M$. By \cite{haefliger62}, the smooth isotopy class of $S$ depends only on $[S] \in H_n(M)$. Since $S$ is Lagrangian, it comes with a canonical formal Legendrian structure (more precisely, a formal Legendrian structure for $\{0\} \times S \subset \bR \times M$, as defined in \cite{murphy12}). Given two homologous Lagrangian spheres, we can use a smooth isotopy between them to compare their canonical formal Legendrian structures. If these coincide, the resulting manifolds $E$ are diffeomorphic, compatibly with the homotopy classes of their almost complex structures. In general, the difference between two formal Legendrian structures for a given $n$-sphere is described by an element of $\pi_{n+1}(V_{n,2n+1},U_n)$, where $V_{n,2n+1}$ is the Stiefel manifold. That homotopy group was analyzed in \cite[Lemmas A.5--A.7]{murphy12}, with the following implications for our situation (compare \cite[Theorem A.4]{murphy12}).

{\em Suppose that $n$ is odd.} Then,
\begin{equation} \label{eq:murphy-map-odd}
\pi_{n+1}(V_{n,2n+1},U_n) \iso \bZ \oplus \bZ.
\end{equation}
A formal Legendrian structure for $S$ gives rise to a stable complex trivialization of $TM|S$. Two such trivializations differ by an element of $\pi_n(U_\infty) \iso \bZ$, and this is one component of \eqref{eq:murphy-map-odd}. For the spheres $S_{\delta_{m+1}}$, all such trivializations are compatible with the stable trivialization of $TM$ coming from the embedding $M \subset \bC^{n+1}$ (because $S_{\delta_{m+1}}$ bounds a Lagrangian disc in $\bC^{n+1}$). Hence, that component of \eqref{eq:murphy-map-odd} is zero in our case. A formal Legendrian structure on $S$ also gives rise to a trivialization of the stabilized normal bundle $\nu_S \oplus \bR$. Two trivializations differ by an element of $\pi_n(O_{n+1})$, and the other component of \eqref{eq:murphy-map-odd} is the image of that element in $\pi_n(S^n) \iso \bZ$. In our construction, this integer is determined by the self-intersection number on $H_{n+1}(E) \iso \bZ$, which by \cite[Eq.~ (9.3)]{maydanskiy-seidel09} depends only on the homology class of $S_{\delta_{m+1}}$. It follows that the formal Legendrian structure contains no additional information. This corrects the proof of \cite[Lemma 1.2]{maydanskiy-seidel09}.

{\em Suppose that $n \geq 4$ is even.} Then,
\begin{equation} \label{eq:murphy-map-even}
\pi_{n+1}(V_{n,2n+1},U_n) \iso \bZ/2.
\end{equation}
Consider the simplest case $m = 2$, and Lagrangian spheres $S_{\delta_{\pm}}$ associated to paths as in Figure \ref{fig:fig}. Let's fix a smooth isotopy between them, and use that to compare their formal Legendrian structures, which leads to an element of \eqref{eq:murphy-map-even}. By embedding $M_2$ into $M_m$ in different ways, one sees that for an isotopy of paths in $\bC$ which crosses over an even number of marked points, there is an associated isotopy of spheres in $M_m$ which respects the formal Legendrian structure. This proves the following:

\begin{lemma} \label{th:lemma}
Suppose that $n=2$. Then, any choice of $\delta_{m+1}$ leads to a manifold $E$ which is diffeomorphic to $T^*S^{n+1}$, and this diffeomorphism is compatible with the homotopy classes of almost complex structures. For higher even $n$, the same holds under the following additional assumption:
\begin{itemize} \item[($\ast$)]
$\delta_{m+1}$ can be connected to a ``standard path'' by an isotopy (rel endpoints) which crosses over an even number of marked points in the plane (here, the ``standard paths'' are the $\delta^{k,l}$ from \cite[Section 7]{maydanskiy-seidel09}). 
\end{itemize}
\end{lemma}
\begin{figure}
\begin{centering}
\begin{picture}(0,0)%
\includegraphics{z-erratum-figure.pstex}%
\end{picture}%
\setlength{\unitlength}{3552sp}%
\begingroup\makeatletter\ifx\SetFigFont\undefined%
\gdef\SetFigFont#1#2#3#4#5{%
  \reset@font\fontsize{#1}{#2pt}%
  \fontfamily{#3}\fontseries{#4}\fontshape{#5}%
  \selectfont}%
\fi\endgroup%
\begin{picture}(1596,919)(1153,-1259)
\put(2401,-1186){\makebox(0,0)[lb]{\smash{{\SetFigFont{11}{13.2}{\rmdefault}{\mddefault}{\updefault}{\color[rgb]{0,0,0}$\delta_-$}%
}}}}
\put(2401,-511){\makebox(0,0)[lb]{\smash{{\SetFigFont{11}{13.2}{\rmdefault}{\mddefault}{\updefault}{\color[rgb]{0,0,0}$\delta_+$}%
}}}}
\end{picture}%
\caption{\label{fig:fig}}
\end{centering}
\end{figure}

The following example shows the importance of assumption ($\ast$). Consider the affine variety
\begin{equation} \label{eq:xy2}
X = \{xy^2 + z_1^2 + \cdots + z_n^2 = 1\} \subset \bC^{n+2}.
\end{equation}
As pointed out in \cite[Example 1.5]{seidel-new}, $X$ (with the standard K{\"a}hler form) is one of the manifolds $E$ constructed in \cite{maydanskiy-seidel09}, corresponding to the choice of path as in the right-hand part of \cite[Figure 2]{maydanskiy-seidel09}. Note that $X$ is a double branched cover of $T^*S^{n+1} = \{xw + z_1^2 + \cdots +z_n^2 = 1\}$, under $w = y^2$.  Take the zero-section $S^{n+1} = \{x = \bar{w}, \, z \in \bR^n\} \subset T^*S^{n+1}$, and let $Z$ be its preimage in $X$. Explicit computation shows that $Z$ is an embedded sphere. The inclusion $Z \hookrightarrow X$ is a homotopy equivalence. By the $h$-cobordism theorem, $X$ must be diffeomorphic to the total space of the normal bundle $\nu_Z$. Again by explicit computation, $\nu_Z$ is the pullback of the normal bundle of the zero-section. From now on, assume that $n$ is even. Then $\nu_Z$ is trivial (since it's classified by twice the class of the tangent bundle of a sphere, inside $\mathit{ker}(\pi_n(\mathit{O}_{n+1}) \rightarrow \pi_n(O_{\infty}))$, and that group is either $\bZ/2$ or zero). Suppose that $n \neq 2,6$. Then, $X \iso S^{n+1} \times \bR^{n+1}$  is not even homeomorphic to $T^*S^{n+1}$, by \cite{james-whitehead}. By comparing this with the argument concerning Figure \ref{fig:fig}, one sees that any isotopy from $S_{\delta_-}$ to $S_{\delta_+}$ necessarily yields a nontrivial obstruction element in \eqref{eq:murphy-map-even}.

It remains to consider the case $n = 6$. Then, for any choice of $\delta_{m+1}$, the resulting $E$ will be diffeomorphic to $T^*S^7 \iso S^7 \times \bR^7$ (one shows this using the $h$-cobordism theorem, and the fact that any $7$-dimensional vector bundle over $S^7$ is trivial). However, there are two possible homotopy classes of almost complex structures ($\pi_7(O_{14}/U_7) \iso \pi_7(O_{\infty}/U_{\infty}) \iso \bZ/2$), and it is not clear which one will arise if ($\ast$) is dropped. In particular, we still don't know what element of \eqref{eq:murphy-map-even} appears there.
%

{\em Acknowledment.} We thank Emmy Murphy for helpful suggestions.